\documentclass{amsart}
\usepackage{shortsalch, xy, eufrak}
\xyoption{all}
\title{Moduli of formal $A$-modules under change of $A$.}
\date{July 2016.}
\begin{document}
\maketitle
\begin{abstract}
We develop methods for computing the restriction map from the cohomology of the automorphism group of a height $dn$ formal group law (i.e., the height $dn$ Morava stabilizer group) to the cohomology of the automorphism group of an $A$-height $n$ formal $A$-module, where $A$ is the ring of integers in a degree $d$ field extension of $\mathbb{Q}_p$. We then compute this map for the quadratic extensions of $\mathbb{Q}_p$ and the height $2$ Morava stabilizer group at primes $p>3$. We show that the these automorphism groups of formal modules are closed subgroups of the Morava stabilizer groups, and we use local class field theory to identify the automorphism group of an $A$-height $1$-formal $A$-module with the ramified part of the abelianization of the absolute Galois group of $K$, yielding an action of $\Gal(K^{ab}/K^{nr})$ on the Lubin-Tate/Morava $E$-theory spectrum $E_2$ for each quadratic extension $K/\mathbb{Q}_p$. Finally, we run the associated descent spectral sequence to compute the $V(1)$-homotopy groups of the homotopy fixed-points of this action; one consequence is that, for each element in the $K(2)$-local homotopy groups of $V(1)$, either that element or its dual is detected in the Galois cohomology of the abelian closure of some quadratic extension of $\mathbb{Q}_p$.
\end{abstract}
\tableofcontents

\section{Introduction.}

Let $K$ be a $p$-adic number field with ring of integers $A$, and let $\mathbb{G}_{1/n}^A$ denote an ``$A$-height $n$ formal $A$-module'' over $\overline{\mathbb{F}}_p$, that is, $\mathbb{G}_{1/n}^A$ is a (one-dimensional) formal group law equipped with complex multiplication by $A$, and its underlying formal group law has $p$-height $dn$, where $d = [K : \mathbb{Q}_p]$. (In the base case $K = \mathbb{Q}_p$, 
$\mathbb{G}_{1/n}^{\hat{\mathbb{Z}}_p}$ is simply a $p$-height $n$ formal group law over $\overline{\mathbb{F}}_p$.)
The automorphism group of the underlying formal group law of $\mathbb{G}_{1/n}^A$ is the well-known {\em height $dn$ Morava stabilizer group,} whose cohomology is the input for many spectral sequences computing stable homotopy groups of spheres and other spectra; see~\cite{MR1333942} and~\cite{MR2030586} for one approach, and chapter 6 of~\cite{MR860042} for another. Among the automorphisms of the underlying formal group law $\mathbb{G}_{1/n}^A$, some automorphisms commute with the complex multiplication by $A$, and others do not; hence the automorphism group of $\mathbb{G}_{1/n}^A$ is naturally a subgroup of the height $dn$ Morava stabilizer group. (It is even a {\em closed} subgroup; see Proposition~\ref{construction of map iK} and also this paper's companion and sequel paper,~\cite{cmah7}.)

In this paper we develop methods for computing the restriction map from the cohomology of the height $dn$ Morava stabilizer group to the cohomology of $\Aut(\mathbb{G}_{1/n}^A)$. In particular, in Theorem~\ref{map induced in S(n)} we compute the continuous linear dual Hopf algebra of the group ring of $\Aut(\mathbb{G}_{1/n}^A)$ (recall that, if $A = \hat{\mathbb{Z}}_p$, this linear dual Hopf algebra is called the {\em height $n$ Morava stabilizer algebra}; see section~6.3 of~\cite{MR860042}), as well as the map from the height $dn$ Morava stabilizer algebra to the linear dual of the group ring of $\Aut(\mathbb{G}_{1/n}^A)$, induced by the inclusion of $\Aut(\mathbb{G}_{1/n}^A)$ into the Morava stabilizer group.

The rest of the paper consists of applications of Theorem~\ref{map induced in S(n)}.
We make the cohomology computations for $n=1$ and $d=2$ and $p>3$; that is, for each quadratic extension $K/\mathbb{Q}_p$, we compute the restriction map in cohomology from the cohomology of the height $2$ Morava stabilizer group to the cohomology of $\Aut(\mathbb{G}_{1}^A) \cong A^{\times}$. Here $A$ is again the ring of integers of $K$. Up to isomorphism, there are only quadratic extensions of $\mathbb{Q}_p$, one of which is unramified (and its relevant cohomological computation is Theorem~\ref{unram computation}), and two of which are totally ramified (and their relevant cohomological computations are Theorems~\ref{easy ram computation} and~\ref{harder ram computation}).

By local class field theory, the norm residue symbol map is an isomorphism $A^{\times} \stackrel{\cong}{\longrightarrow} \Gal(K^{ab}/K^{nr})$, where $K^{nr}$ is the compositum of the unramified extensions of $K$, and $K^{ab}$ the compositum of the abelian extensions of $K$; see Theorem~\ref{artin reciprocity} for a quick review of this fact. The natural isomorphism $\Aut(\mathbb{G}_{1}^A) \cong A^{\times}$, composed with the norm residue symbol, embeds $A^{\times}$ as a closed (by Proposition~\ref{construction of map iK}) subgroup of the height $2$ Morava stabilizer group; hence, by the work of Goerss and Hopkins (see~\cite{MR2125040}), there exists an action of $\Gal(K^{ab}/K^{nr})$ on (a model for) the Lubin-Tate/Morava $E$-theory spectrum $E_2$, for each quadratic extension $K/\mathbb{Q}_p$, and by~\cite{MR2030586}, there exists a descent spectral sequence whose input is the (continuous) Galois cohomology of $K^{ab}/K^{nr}$ and whose output is the homotopy-fixed point spectrum $E_2^{h\Gal(K^{ab}/K^{nr})}$.
In Theorem~\ref{topological computation} we run this spectral sequence after smashing with the Smith-Toda complex $V(1)$ at each prime $p>3$, and we compute the resulting map from the homotopy groups of $L_{K(2)}V(1)$ to the homotopy groups of the homotopy fixed-point spectrum $V(1)\smash E_2^{h\Gal(K^{ab}/K^{nr})}$. 

One interesting consequence is that the map 
\begin{equation}\label{unram map 430} \pi_*(L_{K(2)}V(1)) \rightarrow \pi_*\left(V(1)\smash E_2^{h\Gal(\mathbb{Q}_p(\zeta_{p^2-1})^{ab}/\mathbb{Q}_p(\zeta_{p^2-1})^{nr})\rtimes \Gal(\mathbb{F}_{p^2}/\mathbb{F}_p)}\right)\end{equation}
is injective on the sub-$\mathbb{F}_p[v_2^{\pm 1}]$-module of $\pi_*(L_{K(2)}V(1))$ generated by $1$ and the element $\zeta_2$ from Hopkins' chromatic splitting conjecture, and the map~\ref{unram map 430} is zero on the rest of $\pi_*(L_{K(2)}V(1))$. Here $\zeta_{p^2-1}$ denotes a primitive $(p^2-1)$th root of unity.

Another interesting consequence is Corollary~\ref{artin product map cor}: 
the product of the restriction maps 
\begin{equation}\label{artin product map 0} 
H^{*}_c(\Aut(\mathbb{G}_{1/2}^{\hat{\mathbb{Z}}_p}); \mathbb{F}_{p^2})
 \rightarrow \prod_{[K: \mathbb{Q}_p] = 2} H^{*}_c(\Gal(K^{ab}/K^{nr}); \mathbb{F}_{p^2})\end{equation}
is injective in cohomological degrees $\leq 1$, and for each homogeneous element\linebreak $x\in H^{*}_c(\Aut(\mathbb{G}_{1/2}^{\hat{\mathbb{Z}}_p}); \mathbb{F}_{p^2})$,
either $x$ or the Poincar\'{e} dual of $x$ has nonzero image under the map~\ref{artin product map 0}. (The product in the map~\ref{artin product map 0} is taken over all isomorphism classes of quadratic extensions of $\mathbb{Q}_p$.) More generally: given a homogeneous element $x\in \pi_*(L_{K(2)}V(1))$, either $x$ or the Poincar\'{e} dual class of $x$ is detected in $\pi_*(E_2^{\Gal(K^{ab}/K^{nr})})$ for some quadratic extension $K/\mathbb{Q}$. I do not know yet if this phenomenon generalizes to higher heights or to smaller primes.

This paper is a complete (and much improved) rewrite of 
much older material I wrote when I was in graduate school.
I am grateful to T. Lawson for suggesting a Galois descent argument used in the proof of Theorem~\ref{topological computation}, and to D. Ravenel for teaching me a great deal about formal modules and stable homotopy when I was a graduate student.

\begin{conventions}
\begin{itemize}
\item In this paper, all formal groups and formal modules are implicitly assumed to be {\em one-dimensional.}
\item Throughout, we will use Hazewinkel's generators for $BP_*$ (and, more generally, for the classifying ring $V^A$ of $A$-typical formal $A$-modules, where $A$ is a discrete valuation ring).
%\item There are two generating sets for $BP_*$ (and, more generally, for the classifying ring $V^A$ of $A$-typical formal $A$-modules, where $A$ is a discrete valuation ring) in common usage: Hazewinkel's generators and Araki's generators. We will write $v_0^A, v_1^A, v_2^A, \dots$ for Araki's generators and $V_0^A,V_1^A, V_2^A , \dots $ for Hazewinkel's generators. Each of them depends on a choice of uniformizer $\pi = v_0^A = V_0^A$ for the discrete valuation ring $A$; see Theorem~\ref{review thm} for a definition of each set of generators. 

%(It is usually not important which set of generators we use, since $v_n^A \equiv V_n^A$ modulo $\pi$; but for integral formulas it matters.)
\item By a ``$p$-adic number field'' we mean a finite field extension of the $p$-adic rationals $\mathbb{Q}_p$ for some prime $p$.
\item When a ground field $k$ is understood from context, we will write $\Lambda( x_1, \dots ,x_n)$ for the exterior/Grassmann $k$-algebra with generators $x_1, \dots ,x_n$, and $P(x_1, \dots ,x_n)$ for the polynomial $k$-algebra with generators $x_1, \dots ,x_n$.
\item We make use of standard conventions when dealing with Hopf algebroids, as in Appendix~1 of~\cite{MR860042}: we write $\eta_L,\eta_R: A\rightarrow \Gamma$ for the left and right unit maps of a Hopf algebroid $(A,\Gamma)$, and if $a\in A$, we sometimes also write $a$ as shorthand for $\eta_L(a)\in \Gamma$.
\item When $\mathbb{G}$ is an affine group scheme over a field $k$, we write $k[\mathbb{G}]^*$ for the Hopf algebra corepresenting $\mathbb{G}$, and given a $k[\mathbb{G}]^*$-comodule $M$, we write $H^*(\mathbb{G}; M)$ for the group scheme cohomology $\Ext_{k[\mathbb{G}]^*-comod}^*(k, M)$.
\item When $G$ is a profinite group and $M$ a discrete $G$-module, we write $H^*_c(G; M)$ for the usual continuous cohomology of $G$, i.e., 
$H^*_c(G; M) = \colim_N H^*(G/N; M^N)$, where the colimit is taken over all finite index normal subgroups $N$ of $G$.
\end{itemize}
\end{conventions}

\section{Moduli of formal $A$-modules under change of $A$.}

The basic definition is:
\begin{definition}\label{def of formal module}
Let $A$ be a commutative ring, and let $R$ be a commutative $A$-algebra. A {\em formal $A$-module over $R$} is a formal group law $G(X,Y) \in R[[X,Y]]$ together with a ring homomorphism $\rho: A \rightarrow \End(G)$ such that $\rho(a)(X) \equiv aX$ modulo $X^2$. 
\end{definition}
The addition in $\End(G)$ is the formal addition given by $G$, and the multiplication is composition. Chapter~21 of~\cite{MR2987372} is a good reference for formal $A$-modules; the paper~\cite{MR745362} is a faster (but more abbreviated) introduction. Another reference which gives at least an attempt at an introductory account is~\cite{cmah1}. 

The classical results on $p$-height and $p$-typicality (as in~\cite{MR0393050}) were generalized to formal $A$-modules, for $A$ a discrete valuation ring (all but the first claim is proven by M. Hazewinkel in chapter~21 of~\cite{MR2987372}; the first claim is easier, and not directly used in this paper, but a proof can be found in~\cite{cmah2}):
\begin{theorem}\label{review thm}
Let $A$ be a discrete valuation ring of characteristic zero, with finite residue field. Then the classifying Hopf algebroid $(L^A,L^AB)$ of formal $A$-modules admits a retract $(V^A,V^AT)$ with the following properties:
\begin{itemize}
\item The inclusion map
$(V^A,V^AT) \hookrightarrow (L^A,L^AB)$
and the retraction map $(L^A,L^AB) \hookrightarrow (V^A,V^AT)$
are maps of graded Hopf algebroids, and are
mutually homotopy-inverse.
\item If $F$ is a formal $A$-module over a commutative $A$-algebra $R$
and the underlying formal group law of $F$ admits a logarithm $\log_F(X)$, 
then the classifying map $L^A \rightarrow R$ factors through the retraction map
$L^A \rightarrow V^A$ if and only if $\log_F(X) = \sum_{n\geq 1} \alpha_n X^{q^n}$
for some $\alpha_1, \alpha_2, \dots \in R\otimes_{\mathbb{Z}} \mathbb{Q}$,
where $q$ is the cardinality of the residue field of $A$.
\item %$V^A \cong A[v_1^A, v_2^A, \dots ]\cong A[V_1^A, V_2^A, \dots ]$ with $v_n^A$ and $V_n^A$ each in grading degree $2(q^n-1)$,
$V^A \cong A[v_1^A, v_2^A, \dots ]$ with $v_n^A$ in grading degree $2(q^n-1)$,
and $V^AT \cong V^A[t_1^A, t_2^A, \dots ]$ with $t_n^A$ in grading degree $2(q^n-1)$.
\item The generators $\{ v_i^{A}\}$ for $V^A$, called the {\em Hazewinkel generators},
are defined as follows: we fix a uniformizer $\pi$ for $A$, and let $v_0^A = \pi$.
The universal $A$-typical formal
$A$-module has logarithm of the form \begin{equation}\label{log equation 1} %\lim_{h\rightarrow\infty} p^{-h}[p^h](x) =  
 \log (x) = \sum_{i\geq 0} \lambda^A_ix^{p^i},\end{equation} 
and the equation
\begin{equation} \label{Hazewinkel relation} \pi\lambda^A_h = \sum_{i=0}^{h-1}\lambda^A_i( v_{h-i}^{A})^{q^i},\end{equation}
can be solved recursively for elements $v_1^A, v_2^A, \dots \in V^A$; these are the Hazewinkel generators.
%\item The generators $\{ v_i^{A}\}$ for $V^A$, called the {\em Araki generators}, are defined as follows: we fix a uniformizer $\pi$ for $A$, and let $v_0^A = \pi$. The universal $A$-typical formal $A$-module has logarithm of the form~\ref{log equation 1}, and the equation \begin{equation} \label{Araki relation}\pi\lambda^A_h = \sum_{i=0}^h\lambda^A_i(v_{h-i}^{A})^{q^i}\end{equation} can be solved recursively for elements $v_1^A, v_2^A, \dots \in V^A$; these are the Araki generators.
\item We have a formula
\begin{eqnarray} %\nonumber \pi \lambda_h^A  & = & \lambda^A_{h-1} (V_1^A)^{q^{h-1}} + \dots  + \lambda^A_1 (V^A)^q_{h-1} + V_h^A, \\
\label{hazewinkels combinatorial formula} \lambda_h^A & = & \sum_{i_1+ \dots + i_r = h} \pi^{-r} v_{i_1}^A(v_{i_2}^A)^{q^{i_1}}
\dots (v_{i_r}^A)^{q^{i_1+\dots +i_{r-1}}},\end{eqnarray}
where $\pi$ is the uniformizer and $q$ the cardinality of the residue field 
of $A$, and all $i_j$ are positive integers.
\end{itemize}
\end{theorem}

\begin{definition}\label{def of typicality and height}
Let $A$ be a discrete valuation ring of characteristic zero, with uniformizer $\pi$, and with finite residue field. Let $R$ be a commutative $A$-algebra, and let $\mathbb{G}$ be a formal $A$-module over $R$.
\begin{itemize}
\item
We say that $\mathbb{G}$ is {\em $A$-typical} if the classifying map $L^A \rightarrow R$ factors through the retraction $L^A \rightarrow V^A$. 
\item 
If $\mathbb{G}$ is $A$-typical and $n$ is a nonnegative integer, we say that $\mathbb{G}$ has {\em $A$-height $\geq n$} if the classifying map $V^A \rightarrow R$ factors through the quotient map $V^A \rightarrow V^A/(\pi, v_1^A, \dots ,v_{n-1}^A)$. We say that $\mathbb{G}$ has {\em $A$-height $n$} if $\mathbb{G}$ has $A$-height $\geq n$ but not $A$-height $\geq n+1$. If $\mathbb{G}$ has $A$-height $\geq n$ for all $n$, then we say that $\mathbb{G}$ has $A$-height $\infty$.
\item 
The inclusion $V^A\rightarrow L^A$ associates, to each formal $A$-module, an $A$-typical formal $A$-module $\typ(\mathbb{G})$ which is isomorphic to it. If $\mathbb{G}$ is an arbitrary (not necessarily $A$-typical) formal $A$-module, we say that $\mathbb{G}$ has {\em $A$-height $n$} if $\typ(\mathbb{G})$ has $A$-height $n$.
\end{itemize}
\end{definition}

The following is proven in~\cite{MR2987372}:
\begin{prop}\label{height of underlying fgl}
Let $p$ be a prime number.
\begin{itemize}
\item Every formal group law over a commutative $\hat{\mathbb{Z}}_p$-algebra admits the unique structure of a formal $\hat{\mathbb{Z}}_p$-module. Consequently, there is an equivalence of categories between formal group laws over commutative $\hat{\mathbb{Z}}_p$-algebras, and formal $\hat{\mathbb{Z}}_p$-modules. Under this correspondence, a formal $\hat{\mathbb{Z}}_p$-module is $\hat{\mathbb{Z}}_p$-typical if and only if its underlying formal group law is $p$-typical. If $\mathbb{G}$ is a formal $\hat{\mathbb{Z}}_p$-module of $\hat{\mathbb{Z}}_p$-height $n$, then its underlying formal group law has $p$-height $n$.
\item If $L,K$ are finite extensions of $\mathbb{Q}_p$ with rings of integers $B,A$ respectively, if $L/K$ is a field extension of degree $d$, and if $\mathbb{G}$ is a formal $B$-module of $B$-height $n$, then the underlying formal $A$-module of $\mathbb{G}$ has $A$-height $dn$.
\item In particular, if $K/\mathbb{Q}_p$ is a field extension of degree $d$, if $K$ has ring of integers $A$, and if $\mathbb{G}$ is a formal $A$-module of $A$-height $n$, then the underlying formal group law of $\mathbb{G}$ has $p$-height $dn$.
Consequently, {\em the only formal groups which admit complex multiplication by $A$ have underlying formal groups of $p$-height divisible by $d$.}
\end{itemize}
\end{prop}

\begin{definition}
Let $\mathbb{G}$ be an $A$-typical formal $A$-module over a commutative $A$-algebra $R$
given by power series $G(X,Y)\in R[[X,Y]]$ and $\rho(X)\in R[[X]]$.
The {\em strict automorphism group scheme of $\mathbb{G}$}, written $\strictAut(\mathbb{G})$, is the group scheme which sends a commutative $A$-algebra $S$ to the group of strict automorphisms of $\mathbb{G}\otimes_RS$, i.e., the group (under composition) of formal power series $f(X)\in S[[X]]$ such that 
\begin{itemize}
\item $f(X)\equiv X\mod X^2$,
\item $f(G(X,Y)) = G(f(X),f(Y))$, and
\item $f(\rho(X)) = \rho(f(X))$.
\end{itemize}

By the usual functor-of-points argument, the strict automorphism scheme $\strictAut(\mathbb{G})$ of $\mathbb{G}$ is co-represented by the Hopf algebra $R\otimes_{V^A}V^AT\otimes_{V^A} R$, where $R$ is an $V^A$-algebra via the ring map $V^A \rightarrow R$ classifying $\mathbb{G}$, and (as is the usual convention, see e.g. chapter~6 of~\cite{MR860042}) $V^AT$ is a left $V^A$-algebra via the left unit map $\eta_L : V^A \rightarrow V^AT$ and $V^AT$ is a right $V^A$-algebra via the right unit map $\eta_R : V^A \rightarrow V^AT$. (Recall that the left and right unit maps classifying the underlying $A$-typical formal $A$-module of the source and target of the universal strict isomorphism of $A$-typical formal $A$-modules.)
We will write $R[\strictAut(\mathbb{G})]^*$ for the co-representing Hopf algebra of $\strictAut(\mathbb{G})$.
\end{definition}
It is worth being careful about notation: $\strictAut(\mathbb{G})$ is a profinite group scheme but often fails to be proconstant, so it is not always the case that 
$R[\strictAut(\mathbb{G})]^*$ is the continuous $R$-linear dual of the group ring
$R[\strictAut(\mathbb{G})(R)]$, even when $R$ is a field.

In Definition~\ref{G notation} we introduce a new notation which we find very convenient:
\begin{definition}\label{G notation}
Let $K$ be a $p$-adic number field with ring of integers $A$ and residue field $k$, and let $n$ be a positive integer.
\begin{itemize}
\item
We write $\tilde{\mathbb{G}}^A_{1/n}$ for the the formal $A$-module over $k[v_n^A]$ classified by the map $V^A \rightarrow k[v_n^A]$ sending $v_n^A$ to $v_n^A$ and sending $v_i^A$ to zero if $i\neq n$.
\item
Let $k^{\prime}$ be a field extension of $k$, and let $\alpha\in (k^{\prime})^{\times}$.
We write ${}_{\alpha}\mathbb{G}^A_{1/n}$ for the the formal $A$-module over $k^{\prime}$ classified by the map $V^A \rightarrow k^{\prime}$ sending $v_n^A$ to $\alpha$ and sending $v_i^A$ to zero if $i\neq n$.
\end{itemize}
\end{definition}

\begin{prop}\label{identification of Sigma}
Let $K$ be a $p$-adic number field with ring of integers $A$. Let $k$ be the residue field of $A$, let $q$ be the cardinality of $k$ and let $\pi$ be a uniformizer for $A$. %Let $k^{\prime}$ be a field extension of $k$, 
Let $n$ be a positive integer.
Then, as a quotient of $V^AT \cong A[v_1^A, v_2^A, \dots][t_1^A, t_2^A, \dots ]$,
the Hopf algebra representing $\strictAut(\tilde{\mathbb{G}}^A_{1/n})$ is
\[ k[v^A_n][t_1^A,t_2^A,\dots]/( t_i^A(v_n^A)^{q^i}-v_n^A(t_i^A)^{q^n} \ \forall i ).\]
\end{prop}
\begin{proof}
In~\cite{MR745362}, Ravenel proves the formula
\begin{equation}\label{right unit formula 44} \sum{}^F_{i,j\geq 0} \eta_L(v_i^A)(t_j^A)^{q^i} \equiv \sum{}^F_{i,j\geq 0} \eta_R(v_i^A)^{q^j}t_j^A \mod \pi ,\end{equation}
where $\sum{}^F$ is the formal sum, i.e., the sum using the formal group law underlying the universal $A$-typical formal $A$-module (the sum is well-defined because there are only finitely many terms in each grading degree). We are following the usual convention that $v_0^A = \pi$ and $t_0^A = 1$.
As a consequence of equation~\ref{right unit formula 44},
\[ V^A\stackrel{\eta_R}{\longrightarrow} k[ v^A_n ]\otimes_{V^A}V^AT = V^AT/\eta_L( v_0^A, v_1^A, \dots , v_{n-2}^A, v_{n-1}^A, v_{n+1}^A, v_{n+2}^A, \dots )\]
is determined by
\begin{equation}\label{eqn 34095} \sum_{i\geq 0}{}^F t_i^A\eta_R(v_n^A)^{q^i} = \sum_{j\geq 0}{}^F v_n^A(t_j^A)^{q^n}.\end{equation}
Matching gradings, we get $t_i^A\eta_R(v_n^A)^{q^i} = v_n^A(t_i^A)^{q^n}$ in
$k[v_n^A]\otimes_{V^A}V^AT\otimes_{V^A} k[v_n^A]$, 
since there is at most one term in each grading degree on each side of the equation~\ref{eqn 34095}.
This gives us the relation in the statement of the theorem.
\end{proof}

\begin{lemma}\label{kernel of forgetful map}
Let $L/K$ be a finite field extension of degree $d$ and ramification degree $e$, with $K,L$ $p$-adic number fields with rings of integers $A,B$ respectively. Let $\ell$ denote the residue field of $B$. % and let $\ell^{\prime}$ be a field extension of $\ell$. 
Then the ring map 
\begin{equation}\label{map of ts} V^AT \cong V^A[t_1^A, t_2^A, \dots ] \rightarrow V^B[t_1^B, t_2^B, \dots ] \cong V^BT\end{equation} classifying the strict formal $A$-module isomorphism underlying the universal strict formal $B$-module isomorphism sends $t_i^A$ to $t_{ie/d}^B$ if $i$ is divisible by the residue degree $d/e$ of $L/K$, and sends $t_i^A$ to zero if $i$ is not divisible by the residue degree $d/e$.

Furthermore, let $n$ be a positive integer. Then the map
\begin{equation}\label{map 4309845} \kappa: V^AT \rightarrow \ell[\strictAut(\tilde{\mathbb{G}}_{1/n}^B)]^*\end{equation}
classifying the universal strict automorphism of $\tilde{\mathbb{G}}_{1/n}^B$
sends $v_{dn}^A$ to $\frac{\pi_A}{\pi_B^e}(v_n^B)^{\frac{q^{ne}-1}{q^n-1}}$, where %$e$ is the ramification degree of $L/K$, 
$q$ is the cardinality of the residue field of $B$, and $\pi_A,\pi_B$ are uniformizers for $A,B$, respectively.
Furthermore, the kernel of the map~\ref{map 4309845}
contains $\eta_L(v_i^A)$ and $\eta_R(v_i^A)$ for all $i\neq dn$. 
\end{lemma}
\begin{proof}
The claim about the behavior of the map~\ref{map of ts} on the generators $t_1^A, t_2^A, \dots $ is a generalization of Lemma~3.11(a) of~\cite{MR745362}.
Proving this claim requires some explanation of how the generators $t_1^A, t_2^A, \dots $ in $V^AT$ work; see e.g. the proof of Theorem~A.2.1.27(d) in~\cite{MR860042}.
Write $\mathbb{G}_{univ}^B$ for the universal $B$-typical formal $B$-module, and write $\log_{\mathbb{G}_{univ}^B}(X) = \sum_{n\geq 0} \lambda_n^B X^{q^n}$ for its logarithm.
Then $\mathbb{G}_{univ}^B$ is the source of the universal strict isomorphism of $B$-typical formal $B$-modules; write
$\mathbb{G}_{univ,0}^B$ for its target. Then $\mathbb{G}_{univ, 0}^B$ has logarithm $\log_{\mathbb{G}_{univ,0}^B}(X) = \sum_{n\geq 0} \eta_R(\lambda_n^B) X^{q^n}$,
where $\eta_R: V^A \rightarrow V^AT$ is the right unit map. The universal strict isomorphism $f: \mathbb{G}_{univ}^B \rightarrow \mathbb{G}_{univ,0}^B$ has inverse given by the formal sum
\begin{equation}\label{universal iso} f(X) = \sum_{n\geq 0}^{\mathbb{G}_{univ}^B} t_n^B X^{q^n}.\end{equation}

The map~\ref{map of ts}
%\begin{align*} \ell[\strictAut({}_{\alpha}\mathbb{G}_{e/dn}^{A_{nr}})]^* = \ell[t_1^{A_{nr}},t_2^{A_{nr}},\dots]/( t_i^{A_{nr}} \alpha^{q^{ei}-1} - (t_{i}^{A_{nr}})^{q^{en}} \ \forall i ) \\ \rightarrow   \ell[t_1^B,t_2^B,\dots]/( t_i^B \beta^{q^i-1}-(t_i^B)^{q^n} \ \forall i )= \ell[\strictAut({}_{\beta}\mathbb{G}_{1/n}^B)]^* \end{align*}
is determined as follows: since the underlying $A_{nr}$-typical formal $A_{nr}$-module of the universal $B$-typical formal $B$-module is the universal $A_{nr}$-typical formal $A_{nr}$-module, the map $\gamma^{\prime}: V^{A_{nr}} \rightarrow V^B$ classifying the underlying $A_{nr}$-typical formal $A_{nr}$-module of the universal $B$-typical formal $B$-module sends $\lambda_n^{A_{nr}}$ to $\lambda_n^B$ for all $n$. Hence the map of Hopf algebroids $\gamma: (V^{A_{nr}},V^{A_{nr}}T) \rightarrow (V^B,V^BT)$ sends $\eta_R(\lambda_n^{A_{nr}})$ to $\eta_R(\lambda_n^B)$ for all $n$, i.e., solving the equation~\ref{universal iso} yields that $\gamma(t_n^{A_{nr}}) = t_n^B$.

The unramified case is similar: 
the map $V^A \rightarrow V^{A_{nr}}$ sends $\lambda_n^A$ to $\lambda_{ne/d}^{A_{nr}}$ if $n$ is divisible by the residue degree $d/e = [K_{nr}: K]$, and sends
$\lambda_n^A$ to zero if $n$ is not divisible by the residue degree, and solving equation~\ref{universal iso} yields that the map
\begin{align*}  
k(\alpha)[\strictAut({}_{\alpha}\mathbb{G}_{1/dn}^A)]^* & \\ = k(\alpha)[t_1^A,t_2^A,\dots]/( t_i^A \alpha^{q^{ei}-1} - (t_{i}^A)^{q^{en}} \ \forall i ) & \rightarrow 
\ell[t_1^{A_{nr}},t_2^{A_{nr}},\dots]/( t_i^{A_{nr}} \alpha^{q^{ei}-1} - (t_{i}^{A_{nr}})^{q^{en}} \ \forall i ) \\ & \ \ \ \ \ \ \ \ \ \ \ = \ell[\strictAut({}_{\alpha}\mathbb{G}_{e/dn}^{A_{nr}})]^* 
\end{align*}
sends $t_n^A$ to $t_{ne/d}^{A_{nr}}$ if $n$ is divisible by the residue degree $d/e = [K_{nr}: K]$, and sends
$t_n^A$ to zero if $n$ is not divisible by the residue degree.

Now for the claim about the kernel of the map~\ref{map 4309845}:
we need to break the problem into two parts, an unramified part and a totally ramified part.
Let $K_{nr}$ denote the maximal unramified extension of $K$ contained in $L$, and write $A_{nr}$ for the ring of integers of $K_{nr}$.
From Lemma~3.11(b) of~\cite{MR745362} we have that the map $\gamma: V^A \rightarrow V^{A_{nr}}$ classifying the underlying $A$-typical formal $A$-module of the universal $A_{nr}$-typical formal $A_{nr}$-module sends $v_n^A$ to $v_{ne/d}^{A_{nr}}$ if the residue degree $d/e$ divides $n$, and $\gamma(v_n^A) = 0$ if $d/e$ does not divide $n$. 

Meanwhile, 
the map $\gamma^{\prime}: V^{A_{nr}} \rightarrow V^B$ can be computed using equation~\ref{Hazewinkel relation}:
\begin{equation}\label{eq 043284} \gamma^{\prime}\left( \sum_{i=0}^{h-1}\lambda^{A_{nr}}_i(v_{h-i}^{A_{nr}})^{q^i} \right) 
 = \frac{\pi_A}{\pi_B} \sum_{i=0}^{h-1}\lambda^{B}_i(v_{h-i}^{B})^{q^i} \end{equation}
and the fact that $\gamma^{\prime}(\lambda_i^{A_{nr}}) = \lambda_i^B$.
Modulo $v_0^B, v_1^B, \dots ,v_{n-2}^B, v_{n-1}^B, v_{n+1}^B, v_{n+2}^B, \dots $,
equation~\ref{hazewinkels combinatorial formula} reads
\[ \lambda_h^B  =  \left\{ \begin{array}{ll} 
\pi_B^{-h/n} (v_n^B)^{\frac{q^h-1}{q^n-1}} & \mbox{\ \ if\ \ } n\mid h \\
0 & \mbox{\ \ if\ \ } n\nmid h, \end{array}\right. \]
and consequently
equation~\ref{eq 043284} reads:
\begin{align*} 
 \sum_{i=0}^{h-1} \lambda_i^B \gamma^{\prime}(v_{h-i}^{A_{nr}})^{q^i} 
  &= \frac{\pi_A}{\pi_B} \pi_B^{-(h-n)/n} ( v_n^B)^{\frac{q^{h-n}-1}{q^n-1}}(v_n^B)^{q^{h-n}} \\
  &= \frac{\pi_A}{\pi_B^{h/n}} ( v_n^B)^{\frac{q^{h}-1}{q^n-1}}.\end{align*}
Now an easy induction gives us that $\gamma^{\prime}(v_h^A)$ has positive $\pi_B$-adic valuation, and hence is zero in $\ell[\strictAut(\tilde{\mathbb{G}}_{1/n}^B)]^*$, as long as $h< en$, and when $h = en$ we get the formula
$\gamma^{\prime}(v_{en}^{A_{nr}}) = \frac{\pi_A}{\pi_B^e} (v_n^B)^{\frac{q^{en}-1}{q^n-1}}$,
which, combined with the unramified computation above, immediately yields that the map~\ref{map 4309845} sends $v_{dn}^A$ to $\frac{\pi_A}{\pi_B^e}(v_n^B)^{\frac{q^{ne}-1}{q^n-1}}$, as claimed.

Now for a slightly more involved induction:
suppose we have shown that $\kappa(v_{n(e+a)}^{A_{nr}}) = 0$ for $a = 1, \dots ,j-1$.
Then equations~\ref{hazewinkels combinatorial formula} and~\ref{eq 043284} yield the equation
\begin{align}\nonumber \kappa(v_{n(e+j)}^{A_{nr}}) + \lambda_{jn}^B \kappa(v_{ne}^{A_{nr}})^{q^{nj}} &= \frac{\pi_A}{\pi_B^{e+j}} \left( v_n^B\right)^{\frac{q^{n(e+j)}-1}{q^n-1}},\mbox{\ \ \  i.e.,}\\
\label{eq 3085} \kappa(v_{n(e+j)}^{A_{nr}}) &= \left( \frac{\pi_A}{\pi_B^{e+j}} - \left(\frac{\pi_A}{\pi_B^e}\right)^{q^{nj}}\frac{1}{\pi_B^j}\right) \left( v_n^B\right)^{\frac{q^{n(e+j)}-1}{q^n-1}},\end{align}
and the right-hand side of equation~\ref{eq 3085} is zero, since the scalars $\pi_A,\pi_B$ live in the residue field of $B$, which is the finite field with $q$ elements, hence the $q$th power map is the identity.
So $\kappa(v_{n(e+j)}^{A_{nr}}) = 0$ for all $j>0$; an easy computation using equations~\ref{hazewinkels combinatorial formula} and~\ref{eq 043284} also yields that $\kappa(v_{j}^{A_{nr}})=0$ for $j>en$ not divisible by $n$.

An easy consequence of equation~\ref{right unit formula 44} is that $\eta_L(v_i)$ is congruent to $\eta_R(v_i)$ modulo $(\eta_L(v_0^A),\eta_L(v_1^A), \dots ,\eta_L(v_{i-1}^A))$. Hence $\eta_R(v_i)$ is in the kernel of $\kappa$ for all $i<dn$.
Now we carry out an induction to show that $\kappa(\eta_R(v_i^A)) = 0$ for all $i>dn$ as well: suppose that we have already shown that 
$\kappa(\eta_R(v_{dn+j}^A)) =0$ for $j=1, \dots, i-1$.
Then, reducing formula~\ref{right unit formula 44} modulo the kernel of $\kappa$, we have
\[ \sum{}^F_{j\geq 0} \eta_L(v_{dn}^A)(t_j^A)^{q^{en}} = \sum{}^F_{j\geq 0, i\leq dn} \eta_R(v_i^A)^{q^{je/d}} t_j^A\]
(using the fact that the cardinality of the residue field of $A$ is $q^{e/d}$),
and in grading degree $2(q^{dn+i} - 1)$, this equation reads
\begin{equation}\label{equation 09438} \eta_L(v_{dn}^A)(t_i^A)^{q^{en}} = \eta_R(v_{dn}^A)^{q^{ie/d}}t_i^A +^F \eta_R(v_{dn+i}^A).\end{equation}
If $i$ is not divisible by the residue degree $d/e$ of $L/K$, then we have already shown that $\kappa(t_i^A) = 0$, and consequently equation~\ref{equation 09438} implies that $\kappa(\eta_R(v_{dn+i}^A)) = 0$, as desired. So suppose instead that $i$ is divisible by the residue degree $d/e$.
We already know that 
\[ \kappa(\eta_R(v_{dn}^A)) = \kappa(\eta_L(v_{dn}^A)) = \frac{\pi_A}{\pi_B^e}(v_n^B)^{\frac{q^{ne}-1}{q^n-1}},\]
and that 
\begin{equation}\label{equation 34085}t_{ie/d}^B(v_{n}^B)^{q^{ie/d}}= v_{n}^B(t_{ie/d}^B)^{q^{n}}\end{equation} in $\ell[\strictAut(\tilde{\mathbb{G}}_{1/n}^B)]^*$.
Hence:
\begin{align}
\nonumber   \frac{\pi_A}{\pi_B^e}(v_n^B)^{\frac{q^{ne}-1}{q^n-1}} (t_{ie/d}^B)^{q^{en}} &= \kappa\left(\eta_L(v_{dn}^A)(t_i^A)^{q^{en}}\right) \\
\nonumber   &= \kappa\left(\eta_R(v_{dn}^A)^{q^i}t_i^A +^F \eta_R(v_{dn+i}^A) \right) \\
\label{d3f} &= \kappa(\eta_R(v_{dn}^A)^{q^{ie/d}}t_i^A) +^F \kappa(\eta_R(v_{dn+i}^A)) \\
\label{d4f}   &= \left( \frac{\pi_A}{\pi_B^e}(v_n^B)^{\frac{q^{ne}-1}{q^n-1}}\right)^{q^{ie/d}}t_{ie/d}^B  +^F \kappa(\eta_R(v_{dn+i}^A)),
\end{align}
with equation~\ref{d3f} due to the formal group law on $\ell[\strictAut(\tilde{\mathbb{G}}_{1/n}^B)]^*$ being precisely the one classified by $\kappa$.
We have $(\frac{\pi_A}{\pi_B^e})^q = \frac{\pi_A}{\pi_B^e}$, since $\frac{\pi_A}{\pi_B^e}$ is an element of $\ell\cong \mathbb{F}_q$, and repeated use of equation~\ref{equation 34085} then implies that
\[ \frac{\pi_A}{\pi_B^e}(v_n^B)^{\frac{q^{ne}-1}{q^n-1}} (t_{ie/d}^B)^{q^{en}} = \left( \frac{\pi_A}{\pi_B^e}(v_n^B)^{\frac{q^{ne}-1}{q^n-1}}\right)^{q^{ie/d}}t_{ie/d}^B ,\]
and consequently equation~\ref{d4f} implies that $\kappa(\eta_R(v_{dn+i}^A)) = 0$, completing the inductive step.
\end{proof}

Lemma~\ref{kernel of forgetful map} shows that the map
\[ \kappa: V^AT \rightarrow \ell[\strictAut(\tilde{\mathbb{G}}_{1/n}^B)]^*\]
factors through the projection $V^AT \rightarrow k\otimes_{V^A} V^AT \otimes_{V^A} k = k[\strictAut(\tilde{\mathbb{G}}_{1/dn}^A)]^*$, where $k$ is the residue field of $A$. This gives us a well-defined map of Hopf algebras
\[ k[\strictAut(\tilde{\mathbb{G}}_{1/dn}^A)]^* \rightarrow \ell[\strictAut(\tilde{\mathbb{G}}_{1/n}^B)]^*,\]
which we compute in Theorem~\ref{map induced in S(n)}:
\begin{theorem}\label{map induced in S(n)}
Let $L/K$ be a finite field extension of degree $d$, with $K,L$ $p$-adic number fields with rings of integers $A,B$ respectively. Let $k,\ell$ be the residue fields of $A$ and $B$,
let $e$ be the ramification degree of $L/K$, let $q$ be the cardinality of $\ell$, and let $\pi_A,\pi_B$ be uniformizers for $A,B$, respectively. Let $n$ be a positive integer.

Then the underlying formal $A$-module of $\tilde{\mathbb{G}}_{1/n}^B$ is $\tilde{\mathbb{G}}_{1/dn}^A$. Furthermore, if $\ell^{\prime}$ is a field extension of $\ell$ and $\beta\in (\ell^{\prime})^{\times}$, then the underlying formal $A$-module of ${}_{\beta}\mathbb{G}_{1/n}^B$ is ${}_{\alpha }\mathbb{G}_{1/dn}^A$, where 
\[ \alpha = \frac{\pi_A}{\pi_B^e}\beta^{\frac{q^{en}-1}{q^n-1}}.\]

Furthermore, the ring map
\begin{align} \nonumber k(\alpha)[\strictAut({}_{\alpha}\mathbb{G}_{1/dn}^A)]^* = k(\alpha)[t_1^A,t_2^A,\dots]/( t_i^A \alpha^{q^{ei}-1} - (t_{i}^A)^{q^{en}} \ \forall i ) \\
\label{map 4308443}\rightarrow 
 \ell[t_1^B,t_2^B,\dots]/( t_i^B \beta^{q^i-1}-(t_i^B)^{q^n} \ \forall i )
= \ell[\strictAut({}_{\beta}\mathbb{G}_{1/n}^B)]^* \end{align}
classifying the strict formal $A$-module automorphism of ${}_{\alpha}\mathbb{G}_{1/dn}^A)$ underlying the universal strict formal $B$-automorphism of ${}_{\beta}\mathbb{G}_{1/n}^B$ sends $t_i^A$ to $t_{ie/d}^A$ if $i$ is divisible by the residue degree $d/e$ of $L/K$, and sends $t_i^A$ to zero if $i$ is not divisible by the residue degree $d/e$.
\end{theorem}
\begin{proof}
These claims all follow from Proposition~\ref{identification of Sigma} and Lemma~\ref{kernel of forgetful map}.
\end{proof}

\section{The $n=1, d=2$ case.}

Recall the following computation from Theorem~6.3.22 of~\cite{MR860042} (this computation is also carried out again, in detail, in the present paper's companion paper,~\cite{cmah7}):
\begin{theorem}{\bf (Cohomology of the height $2$ Morava stabilizer group at large primes.)}\label{coh of ht 2 morava gp}
Let $k$ be a finite field of characteristic $p>3$.
Then 
\[ H^{*,*}(\strictAut({}_1\mathbb{G}_{1/2}^{\hat{\mathbb{Z}}_p}); k) \cong
  \Lambda(\zeta_2) \otimes_k k\{ 1,h_{10},h_{11},\eta_2h_{10},\eta_2h_{11},\eta_2h_{10}h_{11}\},\]
with bidegrees as follows:
\begin{equation}\label{degree chart 1}
\begin{array}{llllll}
\mbox{Coh.\ class}          & \mbox{Coh.\ degree} & \mbox{Int.\ degree} \\ % & \mbox{Rav.\ degree} & \mbox{Image\ under\ } \sigma \\
1                           & 0                   & 0                 \\%  & 0                   & 1  \\
h_{10}                      & 1                   & 2(p-1)            \\%  & 1                   & h_{11} \\
h_{11}                      & 1                   & 2p(p-1)           \\%  & 1                   & h_{10} \\
\zeta_2                    & 1                   & 0                 \\%  & 2                   & \zeta_2 \\
h_{10}\eta_2                & 2                   & 2(p-1)           \\%   & 3                   & -h_{11}\eta_2 \\
h_{11}\eta_2                & 2                   & 2p(p-1)          \\%   & 3                   & -h_{10}\eta_2 \\
h_{10}\zeta_2                & 2                   & 2(p-1)           \\%   & 3                   & h_{11}\zeta_2 \\
h_{11}\zeta_2                & 2                   & 2p(p-1)          \\%   & 3                   & h_{10}\zeta_2 \\
h_{10}h_{11}\eta_2          & 3                   & 0                 \\%  & 4                   & h_{10}h_{11}\eta_2 \\
h_{10}\eta_2\zeta_2         & 3                   & 2(p-1)           \\%   & 5                   & -h_{11}\eta_2\zeta_2 \\
h_{11}\eta_2\zeta_2         & 3                   & 2p(p-1)          \\%   & 5                   & -h_{10}\eta_2\zeta_2 \\
h_{10}h_{11}\eta_2\zeta_2    & 4                   & 0                \\%   & 6                   & h_{10}h_{11}\eta_2\zeta_2 .
    \end{array} \end{equation}
where the cup products in $\mathbb{F}_p\{ 1,h_{10},h_{11},h_{10}\eta_2,h_{11}\eta_2,h_{10}h_{11}\eta_2 \}$ are all zero aside from the Poincar\'{e} duality cup products, i.e.,
each class has the obvious dual class such that the cup product of the two is
$h_{10}h_{11}\eta_2$, and the remaining cup products are all zero.
The internal/topological degrees are defined modulo $\left| v_2\right| = 2(p^2-1)$.

In the cobar complex for the Hopf algebra $\mathbb{F}_p[{}_{1}\mathbb{G}_{1/2}^{\hat{\mathbb{Z}}_p}]^*$, we have cocycle representatives:
\begin{align}
\label{cocycle rep 1} h_{10} &= \left[ t_1\right], \\
\label{cocycle rep 2} h_{11} &= \left[ t_1^p\right], \\
\label{cocycle rep 3} \zeta_2 &= \left[ t_2 + t_2^p - t_1^{p+1} \right], \\
\label{cocycle rep 4} h_{10}\eta_2 &= \left[ t_1\otimes t_2 - t_1\otimes t_2^p + t_1\otimes t_1^{p+1} +  t_1^2 \otimes t_1^p \right] \\
\label{cocycle rep 5} h_{11}\eta_2 &= \left[ t_1^p\otimes t_2^p - t_1^p\otimes t_2 + t_1^p\otimes t_1^{p+1} +  t_1^{2p} \otimes t_1 \right] .
\end{align}
\end{theorem}

\begin{prop}
Let $K$ be a $p$-adic number field with rings of integers $A$. Let $k$ be the residue field of $A$, let $\pi$ be a uniformizer for $A$, let $k^{\prime}$ be a field extension of $k$, and let $\alpha\in (k^{\prime})^{\times}$.
Then the profinite group scheme ${}_{\alpha}\mathbb{G}^A_1$ is the proconstant group scheme taking the value $1 + \pi A$, the group (under multiplication) of $1$-units in $A$. That is, the Hopf algebra $k^{\prime}[{}_{\alpha}\mathbb{G}^A_1]^*$ is the continuous $k^{\prime}$-linear dual of the topological group ring $k^{\prime}[1+\pi A]$.
\end{prop}
\begin{proof}
It follows from the Barsotti-Tate module generalization of the well-known Dieudonn\'{e}-Manin classification of $p$-divisible groups over algebraically closed fields (see~\cite{MR0157972}; also see~\cite{MR1393439} for a nice treatment of the theory of Barsotti-Tate modules) that 
${}_{\alpha}\mathbb{G}^A_1\otimes_{k^{\prime}} \overline{k}^{\prime} \cong {}_{\beta}\mathbb{G}^A_1\otimes_k \overline{k}^{\prime}$ for all $\alpha,\beta\in k^{\prime}$, where $\overline{k}^{\prime}$ is the algebraic closure of $k^{\prime}$; and that $\strictAut({}_{\alpha}\mathbb{G}^A_1\otimes_{k^{\prime}} \overline{k}^{\prime})$ is the proconstant group scheme with value $1+\pi A$. So we just need to show that $\strictAut({}_{\alpha}\mathbb{G}^A_1)$ is already proconstant.

Theorem~\ref{map induced in S(n)} gives us that
\[ k^{\prime}[\strictAut({}_{\alpha}\mathbb{G}^A_1)]^* \cong k^{\prime}[t_1, t_2, \dots ]/\left( t_i \alpha^{q^i-1} - t_i^{q}\ \ \forall i \right),\]
where $q$ is the cardinality of $k$. Our argument is essentially the same as that of the proof of Theorem~6.2.3 in Ravenel's book~\cite{MR860042}: an affine profinite group scheme $\dots \rightarrow G_2 \rightarrow G_1 \rightarrow G_0$ over a field $k$ is proconstant if and only if the corepresenting Hopf algebroid $k[G_n]^*$ has a $k$-linear basis $\{ y_i\}_{i\in I}$ of idempotent elements such that $y_iy_j = 0$ for all $i\neq j$ (see e.g.~\cite{MR547117}). In the case of $\strictAut({}_{\alpha}\mathbb{G}^A_1)$, it is profinite by virtue of being the limit (over $n$) of the strict automorphism group scheme of a formal $A$-module $n$-bud, and the strict automorphism group scheme of a formal $A$-module $n$-bud is corepresented by the Hopf algebra
$k^{\prime}[t_1, t_2, \dots ,t_m ]/\left( t_i \alpha^{q^i-1} - t_i^{q}, \ i=1, \dots ,m\right)$, where $m$ is the integer floor of $\log_q n$. That Hopf algebra splits, as a $k^{\prime}$-algebra, as the tensor product of copies of $k^{\prime}[t_i]/(t_i\alpha^{q^i-1} - t_i^q)$ for various $i$. 
The map of $k^{\prime}$-algebras 
$k^{\prime}[s]/(s - s^q) \rightarrow k^{\prime}[t_i]/(t_i\alpha^{q^i-1} - t_i^q)$
sending $s$ to $\alpha^{\frac{1-q^i}{q-1}}t_i$ is an isomorphism of $k^{\prime}$-algebras, and $k^{\prime}[s]/(s - s^q)$ admits a $k^{\prime}$-linear basis of idempotents whose pairwise products are all zero, namely, 
$1-s^{q-1}$ and $-\sum_{j=1}^{q-1} (a^is)^j$ for $i=1, \dots ,q-1$ (this basis is taken from Theorem~6.2.3 of~\cite{MR860042}), where $a$ is any generator for $\mathbb{F}_q^{\times} \cong k^{\times} \subseteq (k^{\prime})^{\times}$.
So $\strictAut({}_{\alpha}\mathbb{G}^A_1)$ is indeed ``already'' (i.e., without any need to change base to an algebraic closure) proconstant over $k^{\prime}$. 
\end{proof}

Theorem~\ref{abelian coh 1} is easy and classical, essentially a part of local class field theory:
\begin{theorem}\label{abelian coh 1}
Let $p>3$,
and let $K/\mathbb{Q}_p$ be a quadratic extension with ring of integers $A$.
Let $\pi$ be a uniformizer for $A$ and let $k$ be the residue field of $A$. 
Then the continuous group cohomology of the profinite group $1+\pi A$ of $1$-units in $A$ is
\[ H^*_c(1 + \pi A; k) \cong \Lambda( h_1,h_2),\]
with $h_1,h_2$ in cohomological degree $1$.
\end{theorem}
\begin{proof}
By Proposition~II.5.5 in~\cite{MR1697859}, the $p$-adic exponential and logarithm maps yield an isomorphism of profinite groups between $1+\pi A$ and the group $\pi A$ under addition. As a profinite abelian group, $\pi A \cong A \cong \hat{\mathbb{Z}}_p \times \hat{\mathbb{Z}}_p$, and it is classical that the continuous cohomology $H^*_c(\hat{\mathbb{Z}}_p; k) = \colim_{n\rightarrow\infty} H^*(\mathbb{Z}/p^n\mathbb{Z}; k)$ is the colimit of the sequence of graded abelian groups
\[\xymatrix{
 H^*(\mathbb{Z}/p\mathbb{Z}; k) \ar[r]\ar[d]^{\cong} &  H^*(\mathbb{Z}/p^2\mathbb{Z}; k) \ar[r]\ar[d]^{\cong} & H^*(\mathbb{Z}/p^3\mathbb{Z}; k) \ar[r]\ar[d]^{\cong} & \dots \\
 \Lambda(h)\otimes_k k[b] \ar[r] & \Lambda(h)\otimes_k k[b] \ar[r] & \Lambda(h)\otimes_k k[b] \ar[r] & \dots }\]
where the horizontal maps send $h$ to $h$ and $b$ to $0$, i.e., $H^*_c(\hat{\mathbb{Z}}_p; k)$ is an exterior algebra on one generator.
\end{proof}

Proposition~\ref{list of quadratic extensions} is well-known, and not difficult; see e.g. section~I.6.6 of~\cite{MR1760253}.
\begin{prop}\label{list of quadratic extensions}
Let $p>2$. Then there are, up to isomorphism, exactly three quadratic extensions of $\mathbb{Q}_p$: the unramified extension $\mathbb{Q}_p(\zeta_{p^2-1})$, where $\zeta_{p^2-1}$ is a primitive $(p^2-1)$th root of unity; and 
two totally ramified extensions $\mathbb{Q}_p(\sqrt{p})$ and $\mathbb{Q}_p(\sqrt{ap})$, where $a$ is any integer satisfying $1\leq a < p$ which does not have a square root in $\mathbb{F}_p$. 
\end{prop}

\begin{theorem}\label{unram computation}
Let $p>3$. Then 
\[ H^{*,*}(\strictAut({}_1\mathbb{G}_1^{\hat{\mathbb{Z}}_p[\zeta_{p^2-1}]}); \mathbb{F}_p) \cong \Lambda(h_{20},h_{21}),\]
with $h_{20},h_{21}$ each in cohomological degree $1$ and internal degree $0$, and 
the restriction map
\begin{align*}
 \Lambda(\zeta_2) \otimes_k k\{ 1,h_{10},h_{11},\eta_2h_{10},\eta_2h_{11},\eta_2h_{10}h_{11}\} & \\
  \cong H^{*,*}(\strictAut({}_1\mathbb{G}_{1/2}^{\hat{\mathbb{Z}}_p}); k) 
  & \stackrel{\res}{\longrightarrow} H^{*,*}(\strictAut({}_1\mathbb{G}_1^{\hat{\mathbb{Z}}_p[\zeta_{p^2-1}]}); \mathbb{F}_p) \\
  & \ \ \ \ \ \ \ \ \ \ \ \ \ \ \ \cong \Lambda(h_{20},h_{21}) ,
\end{align*}
induced by the inclusion of the profinite subgroup 
$\strictAut({}_1\mathbb{G}_1^{\hat{\mathbb{Z}}_p[\zeta_{p^2-1}]})$ of
$\strictAut({}_1\mathbb{G}_{1/2}^{\hat{\mathbb{Z}}_p})$,
 is the map of graded $\mathbb{F}_p$-algebras determined by:
\begin{align*}
 \res(\zeta_2) &= h_{20} + h_{21} ,\\
 \res(h_{10}) &= 0 ,\\
 \res(h_{11}) &= 0 ,\\
 \res(h_{10}\eta_2) &= 0,\\
 \res(h_{11}\eta_2) &= 0,\\
 \res(h_{10}h_{11}\eta_2) &= 0.
\end{align*}
\end{theorem}
\begin{proof}
%We need cocycle representatives for cohomology classes in $H^{*,*}(\strictAut({}_1\mathbb{G}_{1/2}^{\hat{\mathbb{Z}}_p}); \mathbb{F}_p)$  in the cobar complex for the Hopf algebra $\mathbb{F}_p[\strictAut({}_1\mathbb{G}_{1/2}^{\hat{\mathbb{Z}}_p})]^* \cong \mathbb{F}_p[t_1, t_2, \dots ]/(t_i^{p^2} - t_i\ \ \mbox{for\ all} i),$ and similarly, for cohomology classes in $H^{*,*}(\strictAut({}_1\mathbb{G}_1^{\hat{\mathbb{Z}}_p[\zeta_{p^2-1}]}); \mathbb{F}_p)$  in the cobar complex for the Hopf algebra $\mathbb{F}_p[\strictAut({}_1\mathbb{G}_1^{\hat{\mathbb{Z}}_p[\zeta_{p^2-1}]})]^* \cong \mathbb{F}_p[t_2, t_4, \dots ]/(t_i^{p^2} - t_i\ \ \mbox{for\ all} i).$ In the cobar complex for the Hopf algebra $\mathbb{F}_p[\strictAut({}_1\mathbb{G}_{1/2}^{\hat{\mathbb{Z}}_p})]^*$, we have cocycle representatives:
%(these formulas are easily extracted from those in section~6.3 of~\cite{MR860042}; the Hopf algebra $\mathbb{F}_p[\strictAut({}_1\mathbb{G}_{1/2}^{\hat{\mathbb{Z}}_p})]^*$ is called $S(2)$ there), and
We will use cocycle representatives~\ref{cocycle rep 1}, \ref{cocycle rep 2}, \ref{cocycle rep 3}, \ref{cocycle rep 4}, and \ref{cocycle rep 5} for cohomology classes in $H^{*,*}(\strictAut({}_1\mathbb{G}_{1/2}^{\hat{\mathbb{Z}}_p}); k)$.
In the cobar complex for $\mathbb{F}_p[\strictAut({}_1\mathbb{G}_1^{\hat{\mathbb{Z}}_p[\zeta_{p^2-1}]})]^*$, one easily computes that $t_2$ and $t_2^p$ are $1$-cocycles
 (since Theorem~\ref{map induced in S(n)} describes\linebreak $\mathbb{F}_p[\strictAut({}_1\mathbb{G}_1^{\hat{\mathbb{Z}}_p[\zeta_{p^2-1}]})]^*$ as a quotient Hopf algebra of $\mathbb{F}_p[\strictAut({}_1\mathbb{G}_{1/2}^{\hat{\mathbb{Z}}_p})]^*$, we can use Ravenel's formulas for the comultiplication on $S(2) \cong \mathbb{F}_p[\strictAut({}_1\mathbb{G}_{1/2}^{\hat{\mathbb{Z}}_p})]^*$ from section~6.3 of~\cite{MR860042}, and simply reduce them modulo $t_1, t_3, t_5, \dots$ to get the comultiplication on $\mathbb{F}_p[\strictAut({}_1\mathbb{G}_1^{\hat{\mathbb{Z}}_p[\zeta_{p^2-1}]})]^*$) which are, modulo coboundaries, $\mathbb{F}_p$-linearly independent; hence
$t_2,t_2^p$ represent two linearly independent classes in $H^1(\strictAut({}_1\mathbb{G}_1^{\hat{\mathbb{Z}}_p[\zeta_{p^2-1}]}); \mathbb{F}_p)$, and so by Theorem~\ref{abelian coh 1}, the cohomology classes of $t_2,t_2^p$ are a minimal set of $\mathbb{F}_p$-algebra generators for\linebreak $H^{*,*}(\strictAut({}_1\mathbb{G}_1^{\hat{\mathbb{Z}}_p[\zeta_{p^2-1}]}); \mathbb{F}_p)$. We write $h_{20},h_{21}$ for the cohomology classes of $t_2,t_2^p$, respectively. 
Applying Theorem~\ref{map induced in S(n)}, the map $\res$ is simply reduction modulo $t_1,t_3,t_5,\dots$ on cocycle representatives; hence $\res(\zeta_2) = h_{20} + h_{21}$ and $\res$ vanishes on all other generators for the ring $H^{*,*}(\strictAut({}_1\mathbb{G}_{1/2}^{\hat{\mathbb{Z}}_p}); k)$.
\end{proof}

\begin{theorem}\label{easy ram computation}
Let $p>3$. Then 
\[ H^{*,*}(\strictAut({}_1\mathbb{G}_1^{\hat{\mathbb{Z}}_p[\sqrt{p}]}); \mathbb{F}_p) \cong \Lambda(h_{10},h_{20}),\]
with $h_{10},h_{20}$ each in cohomological degree $1$, $h_{10}$ in internal degree $2(p-1)$, and $h_{20}$ in internal degree $0$, and 
the restriction map
\begin{align*}
 \Lambda(\zeta_2) \otimes_k k\{ 1,h_{10},h_{11},\eta_2h_{10},\eta_2h_{11},\eta_2h_{10}h_{11}\} & \\
  \cong H^{*,*}(\strictAut({}_1\mathbb{G}_{1/2}^{\hat{\mathbb{Z}}_p}); k) &
  \stackrel{\res}{\longrightarrow} H^{*,*}(\strictAut({}_1\mathbb{G}_1^{\hat{\mathbb{Z}}_p[\sqrt{p}]}); \mathbb{F}_p) \\
  & \ \ \ \ \ \ \ \ \ \ \ \ \ \ \ \cong \Lambda(h_{10},h_{20}) ,
\end{align*}
induced by the inclusion of the profinite subgroup 
$\strictAut({}_1\mathbb{G}_1^{\hat{\mathbb{Z}}_p[\sqrt{p}]})$ of
$\strictAut({}_1\mathbb{G}_{1/2}^{\hat{\mathbb{Z}}_p})$,
 is the map of graded $\mathbb{F}_p$-algebras determined by:
\begin{align*}
 \res(\zeta_2) &=  2h_{20},\\
 \res(h_{10}) &= h_{10} ,\\
 \res(h_{11}) &= h_{10} ,\\
 \res(h_{10}\eta_2) &= 0,\\
 \res(h_{11}\eta_2) &= 0,\\
 \res(h_{10}h_{11}\eta_2) &= 0.
\end{align*}
\end{theorem}
\begin{proof}
Very similar computation to Theorem~\ref{unram computation}. 
We use cocycle representatives~\ref{cocycle rep 1}, \ref{cocycle rep 2}, \ref{cocycle rep 3}, \ref{cocycle rep 4}, and \ref{cocycle rep 5} for the $\mathbb{F}_p$-algebra generators of $H^{*,*}(\strictAut({}_1\mathbb{G}_{1/2}^{\hat{\mathbb{Z}}_p}); k)$.
Theorem~\ref{map induced in S(n)} tells us that, as a quotient of the Hopf algebra \[ \mathbb{F}_p[\strictAut({}_1\mathbb{G}_{1/2}^{\hat{\mathbb{Z}}_p})]^* \cong \mathbb{F}_p[t_1, t_2, \dots ]/(t_i^{p^2} - t_i),\] the Hopf algebra 
$\mathbb{F}_p[\strictAut({}_1\mathbb{G}_1^{\hat{\mathbb{Z}}_p[\sqrt{p}]})]^*$ is $\mathbb{F}_p[t_1, t_2, \dots ]/(t_i^p - t_i)$.
The $1$-cocycles $t_1$ and $t_2 - \frac{1}{2}t_1^2$ in the cobar complex for 
$\mathbb{F}_p[\strictAut({}_1\mathbb{G}_1^{\hat{\mathbb{Z}}_p[\sqrt{p}]})]^*$
are easily seen to be $\mathbb{F}_p$-linearly independent modulo coboundaries, so Theorem~\ref{abelian coh 1} again implies that the cohomology classes of 
$t_1$ and $t_2 - \frac{1}{2}t_1^2$ are a minimal set of $\mathbb{F}_p$-algebra generators for \linebreak $H^{*,*}(\strictAut({}_1\mathbb{G}_1^{\hat{\mathbb{Z}}_p[\sqrt{p}]}); \mathbb{F}_p)$.

Reducing cocycles~\ref{cocycle rep 1}, \ref{cocycle rep 2}, and \ref{cocycle rep 3} modulo $(t_1^p - t_1,t_2^p-t_2)$ immediately yields the given formulas for $\res(\zeta_2), \res(h_{10})$, and $\res(h_{11})$. For $\res(h_{10}\eta_2)$, we see that reducing~\ref{cocycle rep 4} modulo $(t_1^p - t_1,t_2^p-t_2)$ yields the $2$-cocycle $t_1\otimes t_1^2 + t_1^2\otimes t_1$, which is the coboundary of $\frac{1}{3} t_1^3$, so $\res(h_{10}\eta_2) = 0$, and a similar computation yields $\res(h_{11}\eta_2) = 0$.
\end{proof}

\begin{theorem}\label{harder ram computation}
Let $p>3$, and 
choose an integer $a$ such that $0<a<p$ and such that $a$ does not have a square root in $\mathbb{F}_p$. 
Then $\mathbb{F}_{p^2}$ has a $(p+1)$st root $\omega$ of $a$, and
the underlying formal $\hat{\mathbb{Z}}_p$-module of ${}_{\omega}\mathbb{G}_1^{\hat{\mathbb{Z}}_p(\sqrt{ap})}$ is ${}_{1}\mathbb{G}_{1/2}^{\hat{\mathbb{Z}}_p}$, and
\[ H^{*,*}(\strictAut({}_{\omega}\mathbb{G}_1^{\hat{\mathbb{Z}}_p(\sqrt{ap})}); \mathbb{F}_p) \cong \Lambda(h_{10},h_{20}),\]
with $h_{10},h_{20}$ each in cohomological degree $1$, $h_{10}$ in internal degree $2(p-1)$, and $h_{20}$ in internal degree $0$, and 
the restriction map
\begin{align*}
 \Lambda(\zeta_2) \otimes_{\mathbb{F}_{p^2}} \mathbb{F}_{p^2}\{ 1,h_{10},h_{11},\eta_2h_{10},\eta_2h_{11},\eta_2h_{10}h_{11}\} & \\
  \cong H^{*,*}(\strictAut({}_{1}\mathbb{G}_{1/2}^{\hat{\mathbb{Z}}_p}\otimes_{\mathbb{F}_p}\mathbb{F}_{p^2}); \mathbb{F}_{p^2}) & \stackrel{\res}{\longrightarrow} H^{*,*}(\strictAut({}_{\omega}\mathbb{G}_1^{\hat{\mathbb{Z}}_p[\sqrt{ap}]}); \mathbb{F}_{p^2}) \\
 & \ \ \ \ \ \ \ \ \ \ \ \ \ \ \ \cong \Lambda(h_{10},h_{20}) ,
\end{align*}
induced by the inclusion of the profinite subgroup 
$\strictAut({}_{\omega}\mathbb{G}_{1}^{\hat{\mathbb{Z}}_p[\sqrt{p}]})$ of
$\strictAut({}_{1}\mathbb{G}_{1/2}^{\hat{\mathbb{Z}}_p}\otimes_{\mathbb{F}_p}\mathbb{F}_{p^2})$,
 is the map of graded $\mathbb{F}_p$-algebras determined by:
\begin{align*}
 \res(\zeta_2) &=  2h_{20},\\
 \res(h_{10}) &= h_{10} ,\\
 \res(h_{11}) &= \omega^{p-1}h_{10} ,\\
 \res(h_{10}\eta_2) &= 0,\\
 \res(h_{11}\eta_2) &= 0,\\
 \res(h_{10}h_{11}\eta_2) &= 0.
\end{align*}
\end{theorem}
\begin{proof}
Very similar computation to Theorem~\ref{easy ram computation}. 
The existence of $\omega$ in $\mathbb{F}_{p^2}$ is very easy: $a$ is a $(p-1)$st root of unity since $a\in \mathbb{F}_p$, hence $\omega$ is a $(p+1)(p-1)$st root of unity, hence $\omega$ is fixed by the square of the Frobenius on the algebraic closure $\overline{\mathbb{F}}_p$. So $\omega\in \mathbb{F}_{p^2}$.

We use cocycle representatives~\ref{cocycle rep 1}, \ref{cocycle rep 2}, \ref{cocycle rep 3}, \ref{cocycle rep 4}, and \ref{cocycle rep 5} for the $\mathbb{F}_{p^2}$-algebra generators of $H^{*,*}(\strictAut({}_{1}\mathbb{G}_{1/2}^{\hat{\mathbb{Z}}_p}\otimes_{\mathbb{F}_p}\mathbb{F}_{p^2}); \mathbb{F}_{p^2})$.
Theorem~\ref{map induced in S(n)} tells us that, as a quotient of the Hopf algebra $\mathbb{F}_p[\strictAut({}_{1}\mathbb{G}_{1/2}^{\hat{\mathbb{Z}}_p}\otimes_{\mathbb{F}_p}\mathbb{F}_{p^2})]^* \cong \mathbb{F}_{p^2}[t_1, t_2, \dots ]/(t_i^{p^2} - t_i)$, the Hopf algebra 
$\mathbb{F}_{p^2}[\strictAut({}_{\omega}\mathbb{G}_1^{\hat{\mathbb{Z}}_p[\sqrt{ap}]})]^*$ is $\mathbb{F}_{p^2}[t_1, t_2, \dots ]/(t_i^p - \omega^{p^i-1} t_i)$, i.e.,
\[ \mathbb{F}_{p^2}[t_1, t_2, \dots ]/\left(t_i^p - \omega^{p-1}t_i \mbox{\ for\ } i\mbox{\ odd,\ } t_i^p -t_i\mbox{\ for\ } i\mbox{\ even}\right).\]
The $1$-cocycles $t_1$ and $t_2 - \frac{\omega^{p-1}}{2}t_1^2$ in the cobar complex for 
$\mathbb{F}_{p^2}[\strictAut({}_{\omega}\mathbb{G}_1^{\hat{\mathbb{Z}}_p[\sqrt{ap}]})]^*$
are again easily seen to be $\mathbb{F}_p$-linearly independent modulo coboundaries, so Theorem~\ref{abelian coh 1} again implies that the cohomology classes of 
$t_1$ and $t_2 - \frac{\omega^{p-1}}{2}t_1^2$ are a minimal set of $\mathbb{F}_p$-algebra generators for $H^{*,*}(\strictAut({}_{\omega}\mathbb{G}_1^{\hat{\mathbb{Z}}_p[\sqrt{ap}]}); \mathbb{F}_{p^2})$.

Reducing cocycles~\ref{cocycle rep 1}, \ref{cocycle rep 2}, and \ref{cocycle rep 3} modulo $(t_1^p - \omega^{p-1}t_1,t_2^p-t_2)$ immediately yields the given formulas for $\res(h_{10})$, $\res(h_{11})$, and $\res(\zeta_2)$. 
For $\res(h_{10}\eta_2)$, we see that reducing~\ref{cocycle rep 4} modulo $(t_1^p - \omega^{p-1} t_1,t_2^p-t_2)$ yields the $2$-cocycle $\omega^{p-1}(t_1\otimes t_1^2 + t_1^2\otimes t_1)$, which is the coboundary of $\frac{\omega^{p-1}}{3} t_1^3$, so $\res(h_{10}\eta_2) = 0$, and a similar computation yields $\res(h_{11}\eta_2) = 0$.
\end{proof}

\section{Relations with local class field theory.}

The profinite group scheme $\Aut({}_1\mathbb{G}_{1/n}^{\hat{\mathbb{Z}}_p}\otimes_{\mathbb{F}_p}\mathbb{F}_{p^n})$ is isomorphic to the constant profinite group scheme on the group of units $\mathcal{O}^{\times}_{D_{1/n,\mathbb{Q}_p}}$ in the maximal order in the central division algebra over $\mathbb{Q}_p$ of Hasse invariant $1/n$; see Remark~\ref{brauer remark}. This group of units $\mathcal{O}^{\times}_{D_{1/n,\mathbb{Q}_p}}$ plays an important role in attempts to generalize classical local class field theory to a ``nonabelian local class field theory'' capable of describing the representations of $\Gal(\overline{\mathbb{Q}}_p/\mathbb{Q}_p)$ of degree $>1$ (the degree $1$ representations are already described by classical local class field theory; see Theorem~\ref{artin reciprocity} for the reason why). The book~\cite{MR1876802} is a standard reference for the nonabelian generalizations, and there one can read about the successes that have been had in producing certain nonabelian generalizations of local class field theory which describe $\ell$-adic representations of $\Gal(\overline{\mathbb{Q}}_p/\mathbb{Q}_p)$, for $\ell\neq p$.

It is a much more difficult problem, however, to produce a nonabelian local class field theory which describes $p$-adic representations of $\Gal(\overline{\mathbb{Q}}_p/\mathbb{Q}_p)$, and much less is known about this case; this is one of the major goals of research in $p$-adic Hodge theory. One issue that makes $p$-adic representations more difficult is the dramatic failure of semisimplicity for integral $p$-adic representations of $\Gal(\overline{\mathbb{Q}}_p/\mathbb{Q}_p)$ and related categories of representations: the group $\mathcal{O}^{\times}_{D_{1/n,\mathbb{Q}_p}}$ (whose $\ell$-adic representations are closely related to $\ell$-adic degree $n$ representations of $\Gal(\overline{\mathbb{Q}}_p/\mathbb{Q}_p)$, for $\ell\neq p$, by the local Langlands correspondences; see~\cite{MR1876802}) has a great deal of nonvanishing mod $p$ cohomology, i.e., the cohomology of the height $n$ Morava stabilizer group. 

In this section, motivated by a remark of E. Artin\footnote{``My own belief is that we know it already, though no one will believe me---that whatever can be said about non-Abelian class field theory follows from what we know
now, since it depends on the behavior of the broad field over the intermediate fields---and
there are sufficiently many Abelian cases,'' E. Artin, 1946; the quotation appears in the paper~\cite{cogdellsurvey}, which refers the reader to page~312 of~\cite{durenbook}.}, we compute the restriction map induced in cohomology by the inclusion of $\Gal(K^{ab}/K^{nr})$ as a subgroup of $\mathcal{O}^{\times}_{D_{1/2,\mathbb{Q}_p}}$ for all quadratic extensions $K/\mathbb{Q}_p$, with $p>3$; this inclusion is constructed using local class field theory in Proposition~\ref{construction of map iK}.

Recall (from e.g. Serre's chapter on local class field theory in~\cite{MR0215665}) Artin reciprocity for $p$-adic number fields:
\begin{theorem}\label{artin reciprocity}
Let $K$ be a $p$-adic number field with ring of integers $A$. Let $\pi$ denote a uniformizer for $A$, let $k$ denote the residue field of $A$, and let $K^{ab},K^{tr},K^{nr}$ denote the compositum of all the finite abelian, tamely ramified, and unramified Galois extensions of $K$, respectively, in some fixed algebraic closure for $K$.
Then there exists commutative diagrams of profinite groups with exact rows:
\[\xymatrix{
 1 \ar[r] \ar[d] & A^{\times} \ar[r]\ar[d]_{\cong} & K^{\times} \ar[r]\ar[d]_{\theta} & \mathbb{Z} \ar[d]^{\eta} \ar[r] & 1 \ar[d] \\
 1 \ar[r]        & \Gal(K^{ab}/K^{nr}) \ar[r]             & \Gal(K^{ab}/K) \ar[r]   & \Gal(K^{nr}/K) \ar[r]     & 1  \\
 1 \ar[r] \ar[d] & 1+ \pi A \ar[r]\ar[d]_{\cong} & A^{\times} \ar[r]\ar[d]_{\cong} & k^{\times} \ar[d]_{\cong} \ar[r] & 1 \ar[d] 
 \\
 1 \ar[r]        & \Gal(K^{ab}/K^{tr}) \ar[r]             & \Gal(K^{ab}/K^{nr}) \ar[r]   & \Gal(K^{tr}/K^{nr}) \ar[r]     & 1 ,
}\]
where $1+\pi$ is the group (under multiplication) of $1$-units in $A$, where $\theta$ is the limit, over all finite abelian extensions $L/K$, of the norm residue symbol maps $K^{\times}/N_{L/K}L^{\times} \stackrel{\cong}{\longrightarrow} \Gal(L/K)$, and where $\eta$ agrees with the embedding of $\mathbb{Z}$ into its profinite completion under the composite
\[ \mathbb{Z} \rightarrow \Gal(K^{nr}/K) \stackrel{\cong}{\longrightarrow} \Gal(\overline{k}/k) \stackrel{\cong}{\longrightarrow} \hat{\mathbb{Z}}.\]
\end{theorem}

\begin{remark}\label{brauer remark}
By the computation of the Brauer group $\Br(K) \cong H^2_c(\Gal(K^{sep}/K); (K^{sep})^{\times})\cong \mathbb{Q}/\mathbb{Z}$ of a local field $K$ (again, see Serre's chapter on local class field theory in~\cite{MR0215665}), every finite-rank central division algebra over $K$ is classified, up to isomorphism, by some rational number $r/s\in \mathbb{Q}/\mathbb{Z}$, called the {\em Hasse invariant} of the division algebra; if $r/s$ is a reduced fraction, then $s^2$ is the rank of the division algebra. It is well-known (again, see Serre's chapter on local class field theory in~\cite{MR0215665}, or~\cite{MR0157972}) that the endomorphism ring of an $A$-height $n$ formal $A$-module over the algebraic closure of $k$ is isomorphic to the maximal order (i.e., maximal compact subring) in the central division algebra with Hasse invariant $1/n$; and furthermore every degree $n$ field extension of $K$ embeds, by a ring homomorphism, into that division algebra. Here $A$ is the ring of integers of $K$ and $k$ is the residue field of $A$.

Since the profinite group scheme $\Aut({}_1\mathbb{G}_{1/n}^{\hat{\mathbb{Z}}_p})$ is already proconstant after base change to $\mathbb{F}_{p^n}$ (this is Theorem~6.2.3 of~\cite{MR860042}), it is unnecessary to base change all the way to $\overline{\mathbb{F}}_p$; the above results from local class field theory provide an embedding of $A^{\times} \cong \Gal(K^{ab}/K^{nr})$ into $\Aut({}_1\mathbb{G}_{1/n}^{\hat{\mathbb{Z}}_p}\otimes_{\mathbb{F}_p} \mathbb{F}_{p^n})$ for each field extension $K/\mathbb{Q}_p$ of degree $n$. Here $A$ again denotes the ring of integers of $K$.
\end{remark}

\begin{prop}\label{construction of map iK}
Let $K/\mathbb{Q}_p$ be a field extension of degree $n$. 
Let $A$ denote the ring of integers of $K$, and let $\pi$ denote a uniformizer for $A$ and $k$ the residue field of $A$. Let $q$ be the cardinality of $k$, and let $\omega$ denote a $\frac{q^{en}-1}{q^n-1}$th root of $\frac{\pi^e}{p}$ in $\mathbb{F}_{p^n}$.
Let 
\[ i_K: \Gal(K^{ab}/K^{nr}) \hookrightarrow \Aut({}_1\mathbb{G}_{1/n}^{\hat{\mathbb{Z}}_p}\otimes_{\mathbb{F}_p}\mathbb{F}_{p^n})\]
be the homomorphism of profinite groups defined as the composite of the norm residue symbol $\Gal(K^{ab}/K^{nr}) \stackrel{\cong}{\longrightarrow} A^{\times} \cong \Aut({}_{\omega}\mathbb{G}_1^{A})$ with the natural embedding 
$\Aut({}_{\omega}\mathbb{G}_1^{A}) \hookrightarrow \Aut({}_1\mathbb{G}_{1/n}^{\hat{\mathbb{Z}}_p})$ of $\Aut({}_{\omega}\mathbb{G}_1^{A})$ as the automorphisms of 
the underlying formal $\hat{\mathbb{Z}}_p$-module of ${}_{\omega}\mathbb{G}_1^{A}$ (which is ${}_1\mathbb{G}_{1/n}^{\hat{\mathbb{Z}}_p}$, by Theorem~\ref{map induced in S(n)}) which preserve the complex multiplication by $A$.

Then the image of $i_K$ is a {\em closed} subgroup of $\Aut({}_1\mathbb{G}_{1/n}^{\hat{\mathbb{Z}}_p}\otimes_{\mathbb{F}_p}\mathbb{F}_{p^n})$.
\end{prop}
\begin{proof}
Let $G_a$ denote the automorphism group of the underlying formal $\hat{\mathbb{Z}}_p$-module $a$-bud of ${}_1\mathbb{G}_{1/n}^{\hat{\mathbb{Z}}_p}\otimes_{\mathbb{F}_p}\mathbb{F}_{p^n}$, so that
$\Aut({}_1\mathbb{G}_{1/n}^{\hat{\mathbb{Z}}_p}\otimes_{\mathbb{F}_p}\mathbb{F}_{p^n})$ is, as a profinite group, the limit of the sequence of finite groups $\dots \rightarrow G_3\rightarrow G_2 \rightarrow G_1$. Let $H_a$ denote the subgroup of 
$\Aut({}_1\mathbb{G}_{1/n}^{\hat{\mathbb{Z}}_p}\otimes_{\mathbb{F}_p}\mathbb{F}_{p^n})$ consisting of those automorphisms whose underlying formal $\hat{\mathbb{Z}}_p$-module $a$-bud automorphism commutes with the complex multiplication by $A$, i.e., those
whose underlying formal $\hat{\mathbb{Z}}_p$-module $a$-bud automorphism is an automorphism of the underlying formal $A$-module $a$-bud of ${}_{\omega}\mathbb{G}_1^A$.
The index of $H_a$ in $\Aut({}_1\mathbb{G}_{1/n}^{\hat{\mathbb{Z}}_p}\otimes_{\mathbb{F}_p}\mathbb{F}_{p^n})$ is at most the cardinality of $G_a$, hence is finite. Now we use the theorem of Nikolov-Segal, from~\cite{MR2276769}: every finite-index subgroup of a topologically finitely generated profinite group is an open subgroup.
The group $\Aut({}_1\mathbb{G}_{1/n}^{\hat{\mathbb{Z}}_p}\otimes_{\mathbb{F}_p}\mathbb{F}_{p^n})$ is topologically finitely generated since it has the well-known presentation
\[ \Aut({}_1\mathbb{G}_{1/n}^{\hat{\mathbb{Z}}_p}\otimes_{\mathbb{F}_p}\mathbb{F}_{p^n}) \cong \left( W(\mathbb{F}_{p^n})\langle S\rangle /(S^n = p, x^{\sigma} S = Sx)\right)^{\times},\]
where $W(\mathbb{F}_{p^n})\langle S\rangle$ is the Witt ring of the field $\mathbb{F}_{p^n}$ with a noncommuting polynomial indeterminate $S$ adjoined, and the relation $x^{\sigma} S = Sx$ is to hold for all $x\in W(\mathbb{F}_{p^n})$, where $x^{\sigma}$ is a lift of the Frobenius automorphism of $\mathbb{F}_{p^n}$ to $W(\mathbb{F}_{p^n})$, applied to $x$; since \[ W(\mathbb{F}_{p^n})\langle S\rangle /\left(S^n = p, x^{\sigma} S = Sx)\right)\] reduced modulo $S$ is a finite ring, $\Aut({}_1\mathbb{G}_{1/n}^{\hat{\mathbb{Z}}_p}\otimes_{\mathbb{F}_p}\mathbb{F}_{p^n})$ has a dense finitely generated subgroup.

So $H_a$ is an open subgroup of $\Aut({}_1\mathbb{G}_{1/n}^{\hat{\mathbb{Z}}_p}\otimes_{\mathbb{F}_p}\mathbb{F}_{p^n})$. Every open subgroup of a profinite group is also closed; consequently each $H_a$ is a closed subgroup of $\Aut({}_1\mathbb{G}_{1/n}^{\hat{\mathbb{Z}}_p}\otimes_{\mathbb{F}_p}\mathbb{F}_{p^n})$, and consequently the intersection $\cap_a H_a$ is a closed subgroup of $\Aut({}_1\mathbb{G}_{1/n}^{\hat{\mathbb{Z}}_p}\otimes_{\mathbb{F}_p}\mathbb{F}_{p^n})$. But $\cap_a H_a$ is the group of all formal power series which are automorphisms of ${}_1\mathbb{G}_{1/n}^{\hat{\mathbb{Z}}_p}\otimes_{\mathbb{F}_p}\mathbb{F}_{p^n}$ and whose polynomial truncations, of any length, commute with the complex multiplication by $A$. Consequently $\cap_a H_a = \Aut({}_{\omega}\mathbb{G}_1^A)$ is a closed subgroup of $\Aut({}_1\mathbb{G}_{1/n}^{\hat{\mathbb{Z}}_p}\otimes_{\mathbb{F}_p}\mathbb{F}_{p^n})$.
\end{proof}

\begin{definition}\label{def of twists}
Let $p$ be a prime number, $n$ a positive integer, and $i$ an integer.
The profinite group scheme $\Aut({}_1\mathbb{G}^{\hat{\mathbb{Z}}_p}_{1/n}\otimes_{\mathbb{F}_p} \mathbb{F}_{p^n})$ is proconstant (see Remark~\ref{brauer remark}); in this definition, and in the rest of the paper, we will write $\Aut({}_1\mathbb{G}^{\hat{\mathbb{Z}}_p}_{1/n}\otimes_{\mathbb{F}_p} \mathbb{F}_{p^n})$ to mean the profinite group given by evaluating that group scheme on $\mathbb{F}_{p^n}$ (i.e., $\Aut({}_1\mathbb{G}^{\hat{\mathbb{Z}}_p}_{1/n}\otimes_{\mathbb{F}_p} \mathbb{F}_{p^n})$ is the usual Morava stabilizer group, as in chapter~6 of~\cite{MR860042}).

We write $\mathbb{F}_{p^n}(i)$ for $\mathbb{F}_{p^n}$ with the action of 
$\Aut({}_1\mathbb{G}^{\hat{\mathbb{Z}}_p}_{1/n}\otimes_{\mathbb{F}_p} \mathbb{F}_{p^n})$
given by the $i$th power of the cyclotomic character, i.e., a given element
$z\in \Aut({}_1\mathbb{G}^{\hat{\mathbb{Z}}_p}_{1/n}\otimes_{\mathbb{F}_p} \mathbb{F}_{p^n})$ acts on $\mathbb{F}_{p^n}(i)$ by multiplication by $\overline{z}^i$, where $\overline{z}$ is the image of $z$ under the reduction map
$\Aut({}_1\mathbb{G}^{\hat{\mathbb{Z}}_p}_{1/n}\otimes_{\mathbb{F}_p} \mathbb{F}_{p^n}) \rightarrow \mathbb{F}_{p^n}^{\times}$, i.e., $\overline{z}$ is the leading term of $z$ as a power series in $\mathbb{F}_{p^n}[[X]]$.

In particular: 
$\strictAut({}_1\mathbb{G}^{\hat{\mathbb{Z}}_p}_{1/n}\otimes_{\mathbb{F}_p} \mathbb{F}_{p^n})$ acts trivially on $\mathbb{F}_{p^n}(i)$ for all $i$,
$\Aut({}_1\mathbb{G}^{\hat{\mathbb{Z}}_p}_{1/n}\otimes_{\mathbb{F}_p} \mathbb{F}_{p^n})$ acts trivially on 
$\mathbb{F}_{p^n}(0)$, and
$\mathbb{F}_{p^n}(i) = \mathbb{F}_{p^n}(i + p^n - 1)$.
\end{definition}

Now we use the computations in Theorems~\ref{unram computation}, \ref{easy ram computation}, and~\ref{harder ram computation} to make the same cohomological computations, but with coefficients twisted by powers of the cyclotomic character, as in Definition~\ref{def of twists}. This cohomology with twisted coefficients is exactly what we need in order to compute the $E_2$-terms of the descent spectral sequences of Theorem~\ref{topological computation}.
\begin{theorem}\label{class field thy thm}
Let $p>3$. Then, for each of the three isomorphism classes of quadratic extensions $K/\mathbb{Q}_p$ (see Proposition~\ref{list of quadratic extensions}), 
the cohomology $H^{*,*}_c(\Gal(K^{ab}/K^{nr}); \mathbb{F}_{p^2}(0))$ is an exterior algebra on two generators in cohomological degree $1$, and:
\begin{itemize}
\item if $K/\mathbb{Q}_p$ is unramified, then
$H^{*,*}_c(\Gal(K^{ab}/K^{nr}); \mathbb{F}_{p^2}(i))$ is isomorphic as a bigraded $\mathbb{F}_{p^2}$-vector space to $H^{*,*}_c(\Gal(K^{ab}/K^{nr}); \mathbb{F}_{p^2}(0))$
if $i$ is divisible by $p^2-1$, and $H^{*,*}_c(\Gal(K^{ab}/K^{nr}); \mathbb{F}_{p^2}(i))\cong 0$ if $i$ is not divisible by $p^2-1$; and
\item if $K/\mathbb{Q}_p$ is totally ramified, then
$H^{*,*}_c(\Gal(K^{ab}/K^{nr}); \mathbb{F}_{p^2}(i))$ is isomorphic as a bigraded $\mathbb{F}_{p^2}$-vector space to $H^{*,*}_c(\Gal(K^{ab}/K^{nr}); \mathbb{F}_{p^2}(0))$
if $i$ is divisible by $p-1$, and $H^{*,*}_c(\Gal(K^{ab}/K^{nr}); \mathbb{F}_{p^2}(i))\cong 0$ if $i$ is not divisible by $p-1$.
\end{itemize}
Meanwhile, 
$H^{*,*}_c(\Aut({}_1\mathbb{G}_{1/2}^{\hat{\mathbb{Z}}_p}); \mathbb{F}_{p^2}(i))$ is isomorphic as a bigraded $\mathbb{F}_{p^2}$-vector space to\linebreak $H^{*,*}_c(\Gal(K^{ab}/K^{nr}); \mathbb{F}_{p^2}(0))$
if $i$ is divisible by $p^2-1$, and $H^{*,*}_c(\Gal(K^{ab}/K^{nr}); \mathbb{F}_{p^2}(i))\cong 0$ if $i$ is not divisible by $p^2-1$.

Finally, the inclusion
\[ i_K: \Gal(K^{ab}/K^{nr}) \hookrightarrow \mathcal{O}^{\times}_{D_{1/2,\mathbb{Q}_p}} \cong \Aut({}_1\mathbb{G}_{1/2}^{\hat{\mathbb{Z}}_p}\otimes_{\mathbb{F}_p}\mathbb{F}_{p^2}) \]
constructed in Proposition~\ref{construction of map iK} induces the following restriction map in cohomology:
\begin{description}
\item[$K/\mathbb{Q}_p$ unramified] 
\begin{align*}
 \Lambda(\zeta_2) \otimes_{\mathbb{F}_{p^2}} \mathbb{F}_{p^2}\{ 1,h_{10},h_{11},\eta_2h_{10},\eta_2h_{11},\eta_2h_{10}h_{11}\} & \\
  \cong H^{*,*}_c(\Aut({}_1\mathbb{G}_{1/2}^{\hat{\mathbb{Z}}_p}); \mathbb{F}_{p^2}(i(p^2-1))) 
  & \stackrel{\res}{\longrightarrow} H^{*,*}_c(\Gal(K^{ab}/K^{nr}); \mathbb{F}_{p^2}(i(p^2-1))) \\
  & \ \ \ \ \ \ \ \ \ \ \ \ \ \ \ \cong \Lambda(h_{20},h_{21}) 
\end{align*}
sends $\zeta_2$ to $h_{20} + h_{21}$ and is zero on all other generators.
\item[$K = \mathbb{Q}_p(\sqrt{p})$]
\begin{align*}
 \Lambda(\zeta_2) \otimes_{\mathbb{F}_{p^2}} \mathbb{F}_{p^2}\{ 1,h_{10},h_{11},\eta_2h_{10},\eta_2h_{11},\eta_2h_{10}h_{11}\} & \\
  \cong H^{*,*}_c(\Aut({}_1\mathbb{G}_{1/2}^{\hat{\mathbb{Z}}_p}); \mathbb{F}_{p^2}(i(p^2-1))) 
  & \stackrel{\res}{\longrightarrow} H^{*,*}_c(\Gal(K^{ab}/K^{nr}); \mathbb{F}_{p^2}(i(p^2-1))) \\
  & \ \ \ \ \ \ \ \ \ \ \ \ \ \ \ \cong \Lambda(h_{10},h_{20}) 
\end{align*}
sends $\zeta_2$ to $2h_{20}$, sends $h_{10}$ and $h_{11}$ to $h_{10}$, and is zero on all other generators.
\item[$K = \mathbb{Q}_p(\sqrt{ap})$ for $a$ a nonsquare in $\hat{\mathbb{Z}}_p^{\times}$] 
\begin{align*}
 \Lambda(\zeta_2) \otimes_{\mathbb{F}_{p^2}} \mathbb{F}_{p^2}\{ 1,h_{10},h_{11},\eta_2h_{10},\eta_2h_{11},\eta_2h_{10}h_{11}\} & \\
  \cong H^{*,*}_c(\Aut({}_1\mathbb{G}_{1/2}^{\hat{\mathbb{Z}}_p}); \mathbb{F}_{p^2}(i(p^2-1))) 
  & \stackrel{\res}{\longrightarrow} H^{*,*}_c(\Gal(K^{ab}/K^{nr}); \mathbb{F}_{p^2}(i(p^2-1))) \\
  & \ \ \ \ \ \ \ \ \ \ \ \ \ \ \ \cong \Lambda(h_{10},h_{20}) 
\end{align*}
sends $\zeta_2$ to $2h_{20}$, sends $h_{10}$ to $h_{10}$ and $h_{11}$ to $\omega^{p-1} h_{10}$, and is zero on all other generators. Here $\omega\in\mathbb{F}_{p^2}$ is a $(p+1)$st root of $a$.
\end{description}

See Theorem~\ref{coh of ht 2 morava gp} for the multiplicative structure of 
$H^{*,*}_c(\strictAut({}_1\mathbb{G}_{1/2}^{\hat{\mathbb{Z}}_p}; \mathbb{F}_{p^2}(0))$ and the cohomological degrees of its generators.
\end{theorem}
\begin{proof}
If $A$ has residue field $k$ and $\mathbb{G}$ is a formal $A$-module over a field extension of $k$, then 
we have the short exact sequence of profinite groups
\begin{equation}\label{profinite grp extension} 1 \rightarrow \strictAut(\mathbb{G}) \rightarrow \Aut(\mathbb{G}) \rightarrow k^{\times} \rightarrow 1,\end{equation}
since the strict automorphisms are simply the automorphisms whose leading coefficient (in $k^{\times}$) is one.
So the quotient map
\[ \Aut({}_a\mathbb{G}_1^A)/\strictAut({}_a\mathbb{G}_1^A)
 \rightarrow \Aut({}_1\mathbb{G}_{1/2}^{\hat{\mathbb{Z}}_p})/\strictAut({}_1\mathbb{G}_{1/2}^{\hat{\mathbb{Z}}_p})\]
is simply the monomorphism $k^{\times} \hookrightarrow \mathbb{F}_{p^2}^{\times}$, where $k$ is the residue field of $k$. Consequently:
\begin{description}
\item[if $K/\mathbb{Q}_p$ is unramified] the action of 
$\Aut({}_a\mathbb{G}_1^A)/\strictAut({}_a\mathbb{G}_1^A)$ on $\mathbb{F}_{p^2}(i)$ is trivial if and only if $i$ is divisible by $p^2-1$, and
\item[if $K/\mathbb{Q}_p$ is totally ramified] the action of 
$\Aut({}_a\mathbb{G}_1^A)/\strictAut({}_a\mathbb{G}_1^A)$ on $\mathbb{F}_{p^2}(i)$ is trivial if and only if $i$ is divisible by $p-1$.
\end{description}

The rest follows easily from Theorems~\ref{unram computation}, \ref{easy ram computation}, and~\ref{harder ram computation}, and the (immediately collapsing) Lyndon-Hochschild-Serre spectral sequence for the extension of profinite groups~\ref{profinite grp extension}.
\end{proof}

\begin{corollary}{\bf (How much of the cohomology of the height $2$ Morava stabilizer group is visible in the cohomology of Galois groups?)}\label{artin product map cor}
The product 
\begin{equation}\label{artin product map} 
H^{*}_c(\Aut({}_1\mathbb{G}_{1/2}^{\hat{\mathbb{Z}}_p}); \mathbb{F}_{p^2}(0))
 \rightarrow \prod_{[K: \mathbb{Q}_p] = 2} H^{*}_c(\Gal(K^{ab}/K^{nr}); \mathbb{F}_{p^2}(0))\end{equation}
of the restriction maps from Theorem~\ref{class field thy thm} is injective in cohomological degrees $\leq 1$. Furthermore, for each homogeneous element
\[ x\in H^{*}_c(\Aut({}_1\mathbb{G}_{1/2}^{\hat{\mathbb{Z}}_p}); \mathbb{F}_{p^2}(0)),\]
either $x$ or the Poincar\'{e} dual of $x$ has nonzero image under the map~\ref{artin product map}.
\end{corollary}

\section{Topological consequences.}

As a consequence of Proposition~\ref{construction of map iK}, if $K/\mathbb{Q}_p$ is a degree $n$ extension, then $\Gal(K^{ab}/K^{nr})$ is a {\em closed} subgroup of the Morava stabilizer group $\Aut({}_1\mathbb{G}_{1/n}^{\hat{\mathbb{Z}}_p}\otimes_{\mathbb{F}_p}\mathbb{F}_{p^n})$. Consequently we can use the machinery of~\cite{MR2030586} to construct and compute the homotopy fixed-point spectrum $E_n^{h\Gal(K^{ab}/K^{nr})}$.

We carry this out in the quadratic case for $p>3$:
\begin{theorem}\label{topological computation}
Let $p>3$.%so that the $p$-primary Smith-Toda complex $V(1)$-exists. 
For each of the three isomorphism classes of quadratic extensions of $\mathbb{Q}_p$ (see Proposition~\ref{list of quadratic extensions}), we compute the $V(1)$-homotopy groups of the homotopy fixed-point spectrum \linebreak $E_2^{h\Gal(K^{ab}/K^{nr})\rtimes \Gal(\mathbb{F}_{p^2}/\mathbb{F}_p)}$:
\begin{description}
\item[$K/\mathbb{Q}_p$ unramified] $\pi_*(V(1) \smash E_2^{h\Gal(K^{ab}/K^{nr})\rtimes \Gal(\mathbb{F}_{p^2}/\mathbb{F}_p)}) \cong \Lambda(h_{20},h_{21})\otimes_{\mathbb{F}_p}\mathbb{F}_p[v_2^{\pm 1}]$,
in homotopy degrees $\left| h_{20}\right| = \left| h_{21}\right| = -1$ and $\left| v_2 \right| = 2(p^2-1)$. The natural map from the homotopy groups of the $K(2)$-local Smith-Toda $V(1)$ is the ring map
\begin{align*} \pi_*(L_{K(2)}V(1)) \cong  \Lambda(\zeta_2) \otimes_{\mathbb{F}_p} {\mathbb{F}_p}\{ 1,h_{10},h_{11},\eta_2h_{10},\eta_2h_{11},\eta_2h_{10}h_{11}\} \otimes_{\mathbb{F}_p} \mathbb{F}_p[v_2^{\pm 1}] \\ \rightarrow \Lambda(h_{20},h_{21})\otimes_{\mathbb{F}_p}\mathbb{F}_p[v_2^{\pm 1}] \cong \pi_*(V(1) \smash E_2^{h\Gal(K^{ab}/K^{nr})\rtimes \Gal(\mathbb{F}_{p^2}/\mathbb{F}_p)})\end{align*}
sending $v_2$ to $v_2$, sending $\zeta_2$ to $h_{20} + h_{21}$, and sending 
$h_{10},h_{11},\eta_2h_{10},\eta_2h_{11},$ and $\eta_2h_{10}h_{11}$ to zero.
\item[$K = \mathbb{Q}_p(\sqrt{p})$] $\pi_*(V(1) \smash E_2^{h\Gal(K^{ab}/K^{nr})\rtimes \Gal(\mathbb{F}_{p^2}/\mathbb{F}_p)}) \cong \Lambda(h_{10},h_{20})\otimes_{\mathbb{F}_p}\mathbb{F}_p[b^{\pm 1}]$,
in homotopy degrees $\left| h_{10}\right| = 2p-3$ and $\left| h_{20}\right| = -1$ and $\left| b \right| = 2(p-1)$. The natural map from the homotopy groups of the $K(2)$-local Smith-Toda $V(1)$ is the ring map
\begin{align*} \pi_*(L_{K(2)}V(1)) \cong \Lambda(\zeta_2) \otimes_{\mathbb{F}_p} {\mathbb{F}_p}\{ 1,h_{10},h_{11},\eta_2h_{10},\eta_2h_{11},\eta_2h_{10}h_{11}\} \otimes_{\mathbb{F}_p} \mathbb{F}_p[v_2^{\pm 1}] \\ \rightarrow \Lambda(h_{10},h_{20})\otimes_{\mathbb{F}_p}\mathbb{F}_p[v_2^{\pm 1}] 
 \cong \pi_*(V(1) \smash E_2^{h\Gal(K^{ab}/K^{nr})\rtimes \Gal(\mathbb{F}_{p^2}/\mathbb{F}_p)})\end{align*}
sending $v_2$ to $b^{p+1}$, sending $\zeta_2$ to $2h_{20}$,
sending $h_{10}$ to $h_{10}$, sending $h_{11}$ to $h_{10}b^{p-1}$
and sending $\eta_2h_{10},\eta_2h_{11},$ and $\eta_2h_{10}h_{11}$ to zero.
\item[$K = \mathbb{Q}_p(\sqrt{ap}),$ with $a$ a nonsquare] $\pi_*(V(1) \smash E_2^{h\Gal(K^{ab}/K^{nr})}) \cong \Lambda(h_{10},h_{20})\otimes_{\mathbb{F}_{p^2}}\mathbb{F}_{p^2}[b^{\pm 1}]$,
in homotopy degrees $\left| h_{10}\right| = 2p-3$ and $\left| h_{20}\right| = -1$ and $\left| b \right| = 2(p-1)$. The natural map from the homotopy groups of the $K(2)$-local Smith-Toda $V(1)$ base-changed to $\mathbb{F}_{p^2}$ is the ring map
\begin{align*} \pi_*\left(V(1)\smash E_2^{h\Aut({}_1\mathbb{G}_{1/2}^{\hat{\mathbb{Z}}_p}\otimes_{\mathbb{F}_p}\mathbb{F}_{p^2})}\right) \cong  \Lambda(\zeta_2) \otimes_{\mathbb{F}_{p^2}} \mathbb{F}_{p^2}\{ 1,h_{10},h_{11},\eta_2h_{10},\eta_2h_{11},\eta_2h_{10}h_{11}\} \otimes_{\mathbb{F}_{p^2}} \mathbb{F}_{p^2}[v_2^{\pm 1}] \\ \rightarrow \Lambda(h_{10},h_{20})\otimes_{\mathbb{F}_{p^2}}\mathbb{F}_{p^2}[v_2^{\pm 1}] \cong \pi_*(V(1) \smash E_2^{h\Gal(K^{ab}/K^{nr})})\end{align*}
sending $v_2$ to $b^{p+1}$, sending $\zeta_2$ to $2h_{20}$,
sending $h_{10}$ to $h_{10}$, sending $h_{11}$ to $\omega^{p-1} h_{10}b^{p-1}$ with $\omega$ a $p-1$st root of $a$,
and sending $\eta_2h_{10},\eta_2h_{11},$ and $\eta_2h_{10}h_{11}$ to zero.
\end{description}
\end{theorem}
\begin{proof}
See~\cite{MR1333942} and ~\cite{MR2030586} for the equivalence
\[ L_{K(n)}S \simeq E_n^{h\Aut({}_1\mathbb{G}_{1/n}^{\hat{\mathbb{Z}}_p}\otimes_{\mathbb{F}_p}\mathbb{F}_{p^n})\rtimes\Gal(\mathbb{F}_{p^n}/\mathbb{F}_p)} .\]
Since $V(1)$ is $E(1)$-acyclic, $L_{K(2)}V(1)$ is weakly equivalent to $L_{E(2)}V(1)$, so $L_{K(2)}V(1) \simeq L_{E(2)}V(1) \simeq V(1) \smash L_{E(2)}S$ since $E(2)$-localization is smashing; see~\cite{MR1192553} for the proof of Ravenel's smashing conjecture.
Since $V(1)$ is finite, $(E_n\smash V(1))^{hG} \simeq E_n^{hG} \smash V(1)$, and now we use the $X = V(1)$ case of the conditionally convergent descent spectral sequence (see e.g.~4.6 of~\cite{MR2645058}, or~\cite{MR2030586})
\begin{align*} 
 E_2^{s,t} \cong H^s_c(G; (E_n)_t(X)) & \Rightarrow \pi_{t-s}((E_n\smash X)^{hG}) \\
  d_r: E_r^{s,t} & \rightarrow E_r^{s+r,t+r-1}.\end{align*}

In the case $n=2$ and $X =V(1)$, we have $(E_2)_*\cong W(\mathbb{F}_{p^2})[[u_1]][u^{\pm 1}]$ with $v_1$ acting by $u_1u^{1-p}$, and consequently $(E_2)_*(V(1)) \cong \mathbb{F}_{p^2}[u^{\pm 1}]$. One needs to know the action of $\Aut({}_1\mathbb{G}_{1/2}^{\hat{\mathbb{Z}}_p}\otimes_{\mathbb{F}_p}\mathbb{F}_{p^2})$ or $\Aut({}_{\omega}\mathbb{G}_1^A) \cong \Gal(K^{ab}/K^{nr})$ on $\mathbb{F}_{p^2}[u^{\pm 1}]$ to compute the $E_2$-term of the spectral sequence; but $\Aut({}_1\mathbb{G}_{1/2}^{\hat{\mathbb{Z}}_p}\otimes_{\mathbb{F}_p}\mathbb{F}_{p^2})$ has the finite-index pro-$p$-subgroup 
$\strictAut({}_1\mathbb{G}_{1/2}^{\hat{\mathbb{Z}}_p}\otimes_{\mathbb{F}_p}\mathbb{F}_{p^2})$, and similarly, $\Aut({}_{\omega}\mathbb{G}_1^A) \cong \Gal(K^{ab}/K^{nr})$ has the finite-index pro-$p$-subgroup $\strictAut({}_{\omega}\mathbb{G}_1^A) \cong \Gal(K^{ab}/K^{tr})$. As a pro-$p$-group admits no nontrivial continuous action on a one-dimensional vector space over a field of characteristic $p$, we only need to know the action of $\Aut({}_1\mathbb{G}_{1/2}^{\hat{\mathbb{Z}}_p}\otimes_{\mathbb{F}_p}\mathbb{F}_{p^2})/\strictAut({}_1\mathbb{G}_{1/2}^{\hat{\mathbb{Z}}_p}\otimes_{\mathbb{F}_p}\mathbb{F}_{p^2})\cong \mathbb{F}_{p^2}^{\times}$ and of 
$\Aut({}_{\omega}\mathbb{G}_1^A)/\strictAut({}_{\omega}\mathbb{G}_1^A) \cong \Gal(K^{tr}/K^{nr}) \cong k^{\times}$ on $\mathbb{F}_{p^2}[u^{\pm 1}]$; i.e., for each $j$,
the $\Aut({}_1\mathbb{G}_{1/2}^{\hat{\mathbb{Z}}_p}\otimes_{\mathbb{F}_p}\mathbb{F}_{p^2})$-module $\mathbb{F}_{p^2}\{u^j\}$ is $\mathbb{F}_{p^2}(i)$ for some $i$, as in Definition~\ref{def of twists}; specifically it is $\mathbb{F}_{p^2}(j)$  (see section~1 of~\cite{MR1333942}). 

Hence Theorem~\ref{class field thy thm} provides the $E_2$-term of the descent spectral sequence for $n=2$ and $X = V(1)$ in each of the four cases $G = \Aut({}_1\mathbb{G}_{1/2}^{\hat{\mathbb{Z}}_p}\otimes_{\mathbb{F}_p}\mathbb{F}_{p^2}),$ $G = \Aut({}_1\mathbb{G}_{1/2}^{\hat{\mathbb{Z}}_p}\otimes_{\mathbb{F}_p}\mathbb{F}_{p^2})\rtimes \Gal(\mathbb{F}_{p^2}/\mathbb{F}_p)$, 
$G = \Gal(K^{ab}/K^{nr})$, and $G = \Gal(K^{ab}/K^{nr})\rtimes \Gal(\mathbb{F}_{p^2}/\mathbb{F}_p)$, and in each case, there is a horizontal vanishing line of finite height already at the $E_2$-page of the spectral sequence (this is computed in Theorem~\ref{class field thy thm}), hence the spectral sequence converges strongly.

The $E_2$-term of the descent spectral sequence, along with the map of $E_2$-terms
induced by the inclusion of the closed subgroup $\Gal(K^{nr}/K^{ab}) \subseteq\Aut({}_1\mathbb{G}_{1/2}^{\hat{\mathbb{Z}}_p}\otimes_{\mathbb{F}_p}\mathbb{F}_{p^2})$, is computed in Theorem~\ref{class field thy thm}. In the case of 
$G = \Aut({}_1\mathbb{G}_{1/2}^{\hat{\mathbb{Z}}_p}\otimes_{\mathbb{F}_p}\mathbb{F}_{p^2})$ and $G = \Gal(K^{ab}/K^{nr})$ for $K = \mathbb{Q}_p(\zeta_{p^2-1})$ and $K = \mathbb{Q}_p(\sqrt{p})$, we computed the cohomology of a $\Gal(\mathbb{F}_{p^2}/\mathbb{F}_p)$-form of the Hopf algebra $\mathbb{F}_{p^2}[G]^*$ in Theorems~\ref{coh of ht 2 morava gp}, \ref{unram computation}, and~\ref{easy ram computation}; since the nonabelian Galois cohomology group $H^1(\Gal(\mathbb{F}_{p^2}/\mathbb{F}_p); GL_n(\mathbb{F}_{p^2}))$ classifying $\Gal(\mathbb{F}_{p^2}/\mathbb{F}_p)$-forms of $n$-dimensional $\mathbb{F}_{p^2}$-vector spaces vanishes (this is a well-known generalization of Hilbert's Theorem 90), the invariants of the $\Gal(\mathbb{F}_{p^2}/\mathbb{F}_p)$-action on $H^*_c(G)$ agree, up to isomorphism of graded $\mathbb{F}_p$-vector spaces, with the results of Theorems~\ref{coh of ht 2 morava gp}, \ref{unram computation}, and~\ref{easy ram computation} (this Galois descent argument was suggested to me by T. Lawson). There is no room for differentials in the descent spectral sequences, so $E_2\cong E_{\infty}$ in each spectral sequence.
\end{proof}

In Theorem~\ref{topological computation}, we have indexed $\pi_*(L_{K(2)}V(1))$
with the homotopy degrees
\[ \begin{array}{llllll}
\mbox{Htpy.\ class}          & \mbox{Degree} \\
1                           & 0                \\
h_{10}                      & 2p-3                 \\
h_{11}                      & 2p^2-2p-1                 \\
\zeta_2                    & -1                 \\
h_{10}\eta_2                & 2p-4                 \\
h_{11}\eta_2                & 2p^2-2p-2                \\
h_{10}\zeta_2                & 2p-4               \\
h_{11}\zeta_2                & 2p^2-2p-2               \\
h_{10}h_{11}\eta_2          & -3                \\
h_{10}\eta_2\zeta_2         & 2p-5                \\
h_{11}\eta_2\zeta_2         & 2p^2-2p-3                \\
h_{10}h_{11}\eta_2\zeta_2    & -4.              \\
    \end{array} \]
It is possible that classes with these names (e.g. $\zeta_2$, which as far as I know was named by M. Hopkins in his work on the ``chromatic splitting conjecture''; I do not know of any reference for this in the literature, however) differ by some power of $v_2$ from the classes with these names used by others in the field; the necessary power of $v_2$ is easily found by comparing grading degrees.

\bibliography{/home/asalch/texmf/tex/salch}{}
\bibliographystyle{plain}
\end{document}